\documentclass[12pt,a4paper]{amsart}
\pdfoutput=1
\usepackage{amsmath,amsfonts,amssymb,amscd,amsthm,comment,cite,float}
\usepackage[shortlabels]{enumitem}

\usepackage{color, soul}
\definecolor{gr}{rgb}{0.7, 1, 0.7}
\definecolor{rr}{rgb}{1, 0.7, 0.7}

\usepackage{latexsym}
\usepackage{graphicx}

\theoremstyle{plain} %definition remark
\newtheorem{theorem}{Theorem}[section]
\newtheorem*{mthm}{Main theorem: Second leading eigenvalue of mixture distribution}

\newtheorem{corollary}[theorem]{Corollary}

\theoremstyle{definition} %definition remark
\newtheorem{definition}[theorem]{Definition}
\theoremstyle{remark} %definition remark
\newtheorem{remark}[theorem]{Remark}

\renewcommand{\mathfrak}{\mathbf}

\newcommand{\ignore}[1]{}

\newcommand{\NN}{\mathbb N}
\newcommand{\RR}{\mathbb R}

\newcommand{\cK}{\mathcal K}

\newcommand{\cX}{\mathcal X}
\newcommand{\cY}{\mathcal Y}

\usepackage{bm}

\title[Stochastic similarity matrices for data clustering]{ On some spectral properties of stochastic similarity matrices for data clustering }
\author{Denis Gaidashev}
\address{Uppsala University, Uppsala, Sweden}
\email{gaidash@math.uu.se}

\author{Ralf Pihlstr\"om}
\address{Stockholm Exergi, Uppsala, Sweden}
\email{ralfpihlstrom@gmail.com}

\author{Martin Ryner}
\address{KTH, Stockholm, Sweden}
\email{martinrr@kth.se}

\subjclass[2010]{}
\date{\today}

\begin{document}
\begin{abstract}
  Clustering in image analysis is a central technique that allows to classify elements of an image. %We describe a simple clustering technique that uses the method of similarity matrices.
  We expand upon recent results in spectral analysis for Gaussian mixture distributions, and in particular, provide conditions for the existence of a spectral gap between the leading and remaining eigenvalues for matrices with entries from a Gaussian mixture with two real univariate components. Furthermore, we describe an algorithm in which a collection of image elements is treated as a dynamical system in which the existence of the mentioned spectral gap results in an efficient clustering.  
\end{abstract}
\maketitle

\tableofcontents

The common feature of modern data analysis is the large amount of data which collectively have a very high dimension.  High dimensionality hinders data processing, and, specifically, data clustering.

However, high dimensionality is quite often caused by data representation, while the number of parameters needed to accurately describe the data is much smaller.

{\it Dimensionality reduction} of large data sets can be achieved by treating a collection of elements in a set as a dynamical system. For example, if a data set is an image, then, initially, each element of the image may be assigned a spacial coordinate in $\RR^m$, where $m$ is typically $2$ or $3$, regardless of its specific shape or size (for example, the location of the center of mass), and a set of quantifiers of its properties, e.g. color, shape, size or granularity. In this way we have a system of, say, $n$ particles.  This collection of points is acted  upon by a matrix whose $(i,j)$-th entries are correlation functions that depend on the distance between particles $i$ and $j$ and a measure of similarity of their properties. Such matrices are usually referred to as \emph{similarity matrices}. One application of the matrix to a $m \cdot n$-dimensional vector of coordinates produces a new state of the dynamical system. Before the next iteration, the entries of the matrix are updated. If the entries of the matrix are constructed appropriately then iterations of this dynamical system reduce distances between the elements with  similar properties. This process of clustering can be also understood in terms of dimensionality reduction: if most of the eigenvalues are small in the absolute value, then the systems quickly converges to the span of the eigenvectors of only few dominant eigenvalues. Such dimensionality reduction has been used, e.g. in \cite{CoiLa} to perform $k$-means clustering on a linear subspace of a smaller dimension. Dimensionality reduction techniques have been extensively studied in \cite{CKLM}, \cite{CLLM1}, \cite{CLLM2}, \cite{PHHV}, \cite{NLCK}, \cite{LAF}, among other works.

%\begin{figure}[h]\centering
%\includegraphics[width=\textwidth]{kmeans.jpg}
%\caption{$k$-means algorithm}
%\label{kmeans}
%\end{figure}

In this paper we analyze the spectrum of some similarity matrices, and address one specific, but nonetheless, important issue: that of the size of the gap between the leading eigenvalue and the rest of the spectrum.

Estimates on the size of the leading eigenvalue of the kernel associated to the similarity matrix of a Gaussian mixture $P=\pi_1 P^1+\pi_2 P^2$ (with probabilities $\pi_1+\pi_2=1$), have been provided in \cite{Shi2}. It is clear from the above discussion that the best scenario for dimensionality reduction is the case when the leading eigenvalue is dominant as compared to the second leading eigenvalue.  We address this case with the our main result of the paper: a theoretical bound on the spectral gap between the first and the second leading eigenvalues in a Gaussian mixture in terms of computable properties of univariate distributions in the mixture, such as the two leading eigenvalues and the leading eigenfunctions (see Section $\ref{section_gap}$ for the precise statement of the result). Dynamical clustering and the spectral gap result are illustrated through simulations in Sections $\ref{sims1}$ and $\ref{sims2}$.

\section{Preliminaries}

We will now introduce the dynamical system that we will study.

Let $(y_i^1,y_i^2)$ denote the coordinate of (the center of mass for) the $i$-th element in the image. Set
$$y = ((y_1^1,y_2^1,\ldots, y_n^1),(y_1^2,y_2^2,\ldots,y_n^2)) \in \mathbb{R}^{2n}.$$
Assume that every element carries $l \in \NN$ properties, or parameters: we call 
$$x_i= (x_i^1,x_i^2,\ldots,x_i^l)$$
the parameter of  element $i$. 

The parameter vector of the system will be denoted by $x$:
$$x=(x_1, x_2, \ldots, x_n) \in \RR^{l n}.$$
Then, let $\| \cdot \|$ stand for the $l_2$ metric in $\RR^l$:
%\[
% d_x(x_i,x_j) :=  \sum_{k=1}^l  |x_i^k - x_j^k |^2.
%\]
We will consider the dynamical system  $T_x: \RR^{2 n} \mapsto \RR^{2 n}$,
\begin{align}
\label{iterates} T_x(y)=\cK(x) \cdot y
\end{align}
where $\cK= K \bigotimes K$, and the entries of the matrix $K(x)$ are given by
\begin{align}
	K_{i,j} & = \frac{1}{n} e^{-{\|x_i -x_j\|^2 \over  \omega^2}},
\label{element}
\end{align}
with a fixed $\omega \in \RR_+$. 

We will sometimes refer to the space of parameters $x_i$ as {\it the parameter space} $\cX$, and that of coordinates $y$ as {\it the position space} $\cY$. For the dynamical system that we have described,
$$\cX=\RR^{l}, \quad \cY=\RR^{2n}.$$

 \subsection{Contraction mappings}

 \begin{definition}[Contraction]
 \label{contraction}
 Let $T$ be a continuous map on a complete metric space $X$. We say $T$ is a contraction
 if there is a number $c < 1$ such that
 \[
 	d(T(x),T(y)) \leq c d(x,y)
 \]
 for all $x,y\in X$. 
 \end{definition}
 
 To demonstrate that the dynamical system $T_x$ is a contraction we can use the following fundamental result.
 
 \begin{theorem}(Perron-Frobenius)  \label{perron}
 Let $K$ be a $n \times n$ matrix with $K_{i,j} > 0$
 for $1 \leq i,j \leq n$. There exists a positive real number $r$
 associated to $K$, called the Perron root or the Perron–Frobenius eigenvalue, such that $r$ is an eigenvalue of $K$ and any other eigenvalue  $\lambda$ necessarily satisfies $|\lambda| < r$. 
 Moreover, 
 \[
 	\min_{i}\sum_j K_{i,j} < r < \max_{i} \sum_j K_{i,j}.
 \]
 \end{theorem}
 
 We can now see, immediately, that since the sum
$$
 \max_{i} \sum_{j=1}^n K_{i,j}=  {1 \over n} \max_{i} \sum_{j=1}^n e^{-{\|x_i -x_j\|^2 \over  \omega^2}}  < 1,
$$ 
 as long as not all $x_i=x_j, i \ne j$, the Perron-Frobenius root is strictly less than $1$, and, the linear map $\cK$ is a contraction on $\cY$.

 Successive iterates of an initial condition $y$ under the map $T_x$ converge to a fixed point in $\RR^{2n}$, that is, eventually the positions of all points $y_i$ stabilize.

 This does not imply the existence of cluster, however. {\it The dynamical explanation of clustering is that depending on the distribution law of the parameters $x_i$,  the matrix $\cK$ might have several leading eigenvalues which dominate the rest of the spectrum. This means that the speed of convergence in the complement of the span of the eigenvectors of these leading eigenvalues is very fast, and the iterates of the initial vector $y$ quickly converge to a low-dimensional span of the leading eigenvectors. Visually, this dimensionality reduction exhibits itself as formation of several clusters, that is from a certain point on, the dynamics is described by very few numbers, e.g. coordinates of the centers of mass of those clusters.}

\subsection{Reproducing kernel space}
%https://math.stackexchange.com/questions/1489670/positive-semidefinite-function

\begin{definition}[Kernel]
We define a \emph{kernel} as a symmetric mapping
$$K : \cX \times \cX \to \mathbb{R}.$$
We say the kernel is \emph{positive semi-definite} if $K(x,x) \geq 0$ for all $x \in \cX$.
\end{definition}

We will further consider  the kernel  expressible as an inner product
\[
	K(x_i,x_j) = \langle \phi(x_i), \phi (x_j) \rangle,
\]
where $\phi$ is a nonlinear map from $\cX$ to an inner product space $H$. 

We say that $K$ corresponds \emph{via} $\phi$ to the inner products in $H$. 

We call a Hilbert space $H_K$ a reproducing kernel space, if it consists of real valued functions $f$  defined on $\cX$ where for each $x \in \cX$ the functional $L_x(f) = f(x)$ is bounded in $H_K$.  By the Riesz representation theorem, there exists $K_x \in H_K$  such that
\[	
	L_x(f) = \langle K_x, f\rangle = f(x)\quad
	\forall f \in H_K.
\]
To every reproducing kernel space $H_K$ there corresponds a unique non-negative definite kernel $K$. Conversely, for any non-negative definite kernel $K$ there corresponds a unique Hilbert space that has $K$ as its reproducing kernel. For more details see ~\cite{Aron} and the Moore-Aronszjan theorem.

\subsection{Kernel matrix and kernel operator}

\begin{definition}[Kernel matrix]
Let $K$ be a kernel. We define the associated \emph{kernel matrix} through
\[
	(K_n)_{i,j} = K(x_i,x_j), \ 1 \le i, j \le n 
\]
\end{definition}

\begin{definition}[Kernel operator]
Let $P$ be a probability distribution with density function $p(x)$,
and let $K$ be a kernel function.

We define the \emph{kernel operator}, associated to $K$ as
\[
	K_P f(x) = \int_\cX K(z,x) f(z) p(z) dz
\]
as a mapping from $H_K$ to $H_K$.
\end{definition}
Any eigenfunction $\phi \in H_K$ and the corresponding eigenvalue $\lambda$ for $K_P$ are given through the relation
\[
 \int_\cX K(z,x)\phi (z) p(z) = \lambda \phi (x).
\]

The kernel matrix and operator are related as follows. Let $\lambda_v$ be an eigenvalue, and $v = (v_1,\ldots,v_n)$ an eigenvector of $K_n$. Since
\[
	K_n v = \lambda_v v,
\]
we have that, for each $i = 1,2\ldots, n$
\[
	\frac{1}{n} \sum_{j=1}^n K(x_i,x_j) v_j
	= \frac{\lambda_v}{n} v_i.
\]

If now $x_1,\ldots,x_n$ are samples from a probability distribution with
density $p(x)$, and $v=(\phi(x_1),\ldots,\phi(x_n))$, then
\begin{align}
	\frac{1}{n} \sum_{j=1}^n K(x_i,x_j) v_j
	\approx \int_x K(x,x_i) \phi (x) p(x) dx,
\label{approx}
\end{align}

From (\ref{approx}), it follows that $\lambda_v/n$ approximates the eigenvalue $\lambda$ of the kernel operator with eigenfunction $\phi$. A first formal proof is due to Baker  \cite{BakerC}, showing that $\lambda_v$'s converge to the eigenvalues of the kernel operator as $n\to\infty$. Koltchinskii and Gin\'e \cite{KolGin} refined this result by  showing that, in particular, the convergence rate is $1/\sqrt{n}$ as $n\to\infty$.

%To the kernel $K$ is associated the empirical kernel operator.
%
%\begin{definition}(Empirical kernel operator)
%  Let $x_1,\ldots,x_n$ be an i.i.d. sample from a distribution $P$. The \emph{empirical kernel operator} $K_{P_n}$ is defined through
%  \[
%  K_{P_n} f(x)  = \int K(x,z)f(z) dP_n(z)  = \frac{1}{n} \sum_{i=1}^n K(x,x_i)f(x_i).
%  \]
%\end{definition}
%
%The empirical kernel operator is related (in the  above sense) to the kernel matrix.
%

\subsection{Eigenvalues}

%Spectral Properties of the Kernel Matrix, PhD
%https://pdfs.semanticscholar.org/28a9/4f648bb990d7022a5a83722fc7af0a376a70.pdf
%page 18

%\begin{definition}[Mercer kernel]
%Let $P$ be a probability measure on $\cX$, and $H_\mu(\cX)$ the associated Hilbert space. Given a sequence $(\lambda_i)_{i\in\mathbb{N}}$ with $\lambda_i \geq 0$, and an orthogonal family of unit norm functions $(\phi_i)_{i\in\mathbb{N}}$ with  $\phi_i\in H_\mu(\cX)$, the associated Mercer kernel is
%\[
%K(x,y) = \sum_{i=1}^\infty \lambda_i \phi_i(x)\phi_i(y),
%\]
%\end{definition}
%where we let $T_K f(x) = \int_\cX K(x,z)f(z)dz$
%be the integral operator associated with $K$.
%\end{definition}
Kolthinskii and Gin\'e \cite{KolGin}  provided a way to compare the finitely many eigenvalues of $K_n$ with the infinitely many of the kernel operator $K_P$. We briefly summarize it here.

Let $x_i$, $i \in \NN$, be independently and identically $P$-distributed $\RR$-valued random variables. Set
$$P_n={1 \over n}\sum_{i=1}^n \delta_{x_i}.$$
Let $g \in L^2_{P_n}(\cX,\RR)$ and let $\Omega$ be . The map
$$g \mapsto {1 \over \sqrt{n}} (g(x_1(\omega), \ldots, x_n(\omega))$$  
defines for each $\omega \in \Omega$ an isometry onto a subspace of $\RR$. By means of this isometry $K_{P_n}$ is identified with the following linear operator on $\RR^n$ with the matrix entries
$$K^n_{i,j}={1 \over n} K(x_i,x_j).$$
Furthermore, introduce
$$\tilde{K}^n_{i,j}={1 \over n} \left(  K(x_i,x_j)-\delta^i_j \right).$$

Next, assume that $K_P$ is a Hilbert-Schmidt operator (i.e. $\int_{\cX} K(x,y)^2 d P(x) d P(y) < \infty$). Let the eigenvalues of both $\tilde{K}^n$ and $K_P$ be non-negative and sorted in the non-increasing order, repeated with multiplicity,
\begin{align*}
	\lambda(\tilde{K}^n) &= (l_1,\ldots,l_n), \quad l_1 \geq l_2 \geq \ldots\\
	\lambda(K_P) &= (\lambda_1,\ldots,\lambda_n, \ldots), 	\quad \lambda_1 \geq \lambda_2 \geq \ldots
\end{align*} 
After filling up the first vector with zeros, define the following $l_2$ distance
\[
	\delta_2^2(\lambda(\tilde{K}^n),\lambda(K_P)) = \inf_{\pi\in\sigma(\mathbb{N})} \sum_{i=1}^\infty (l_i-\lambda_{\pi(i)})^2,
\]
where $\sigma(\mathbb{N})$ is the set of all bijections on $\mathbb{N}$. We then have the following.

\begin{theorem}[\cite{KolGin}, Theorem 3.1] Suppose that $K_P$ is a Hilbert-Schmidt kernel operator. Then
   \[
  \delta_2(\lambda(\tilde{K}^n),\lambda(K_P)) \xrightarrow[\text{a.s.}]{} 0.
  \]
\end{theorem}
An important result in the same direction is  the following theorem, 
due to Bonami and Karoui  \cite{BonKar}. 
\begin{theorem}[\cite{BonKar}, Theorem 4]  \label{bound}
  Let $\cX$ be a locally compact metric space, and let $P$ be a probability distribution on $\cX$. Let $K(\cdot,\cdot)$ be a Hermitian kernel, continuous on $\cX \times \cX$, and positive semi-definite. Let $K^n$ be the associated kernel matrix, and let $K_P$ denote the integral operator. If $K(\cdot,\cdot)$
is bounded and if $R := \sup_x |K(x,x)|$ is finite, then we have
\[
\| \lambda(K^n) - \lambda(K_P) \|_{l_2} \leq \frac{ R(\xi + 1) }{\sqrt{n}}
\]
for any $\xi > 0$, with probability at least equal to $1- e^{-\xi ^2}$. 
\end{theorem}

\subsection{Mixture distributions}

\begin{definition}[Mixture distribution]
We call
\[
P = \sum_{g = 1}^G \pi_g P^g
\]
a \emph{mixture distribution} with weights $\pi_g$ and mixing components $P^g$, $g=1,\ldots,G$, and where $\sum_{g=1}^G \pi_g = 1$.

For each mixing component $P^g$, the corresponding operator $K_{P^g}$ is
\[
K_{P^g} f(x) = \int K(x,z)f(z)d P^g (z).
\]
\end{definition}

\begin{figure}[h]\centering
\includegraphics[width=0.6\textwidth]{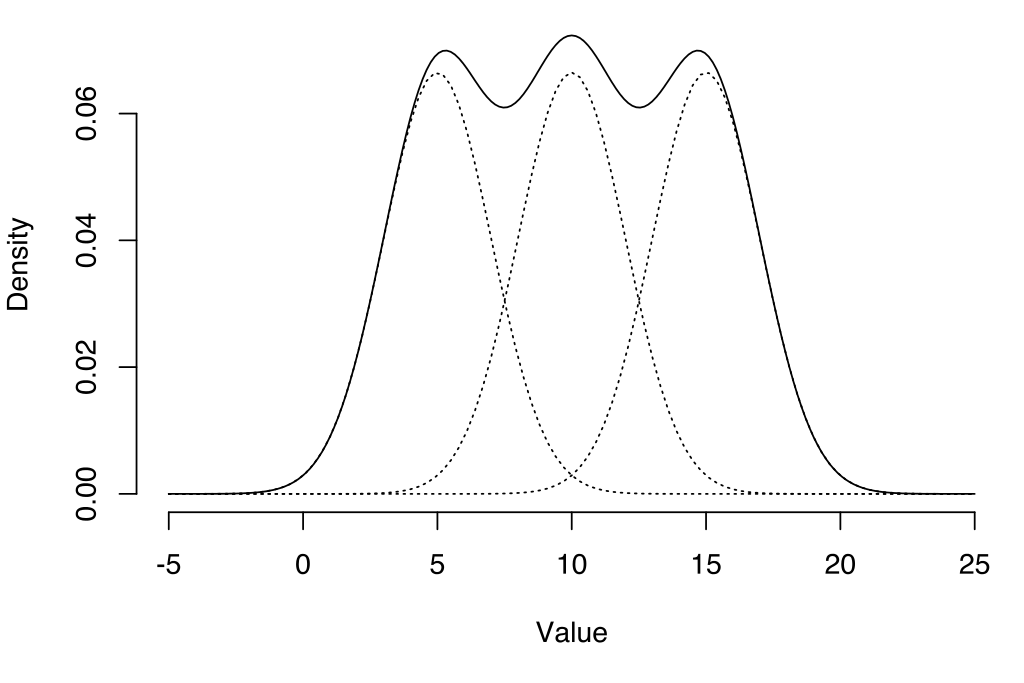}
\caption{Density of a mixture of three normal distributions ($\mu_1 = 5$, 
$\mu_2 = 10$, $\mu_3 = 15$,
$\sigma_1 = \sigma_2 = \sigma_3 = 2$ with equal weights.)}
\label{gme}
\end{figure}

One of the central results about the spectrum of $K_{P^g}$ is an estimate for the top eigenvalue of the kernel operator, due to Shi et al. \cite{Shi2}.

\begin{theorem}[\cite{Shi2}, Theorem 3] \label{shi}
Let $P = \pi_1 P^1 + \pi_2 P^2$ be a mixture distribution on $\mathbb{R}^l$ with $\pi_1 + \pi_2 = 1$. Given a positive semi-definite kernel $K$, denote the top eigenvalues of $K_P$, $K_{p^1}$, $K_{p^2}$ as $\lambda_0$, $\lambda_0^1$, $\lambda_0^2$, respectively.

Then $\lambda_0$ (the top eigenvalue of $K_P$) satisfies
\[
\max ( \pi_1 \lambda_0^1, \pi_2 \lambda_0^2) \leq \lambda_0 \leq  \max ( \pi_1 \lambda_0^1, \pi_2 \lambda_0^2) + r
\]
where
\[
r = \left( \pi_1 \pi_2 \int \int K(x,z)^2 dP^1(x) dP^2(z) \right)^{1/2}.
\]
\end{theorem}

The same authors also provide an estimate on the leading eigenvector of the mixture distribution.

\begin{corollary}[\cite{Shi2}, Corollary 2]\label{cor}
Let $P = \pi_1 P^1 + \pi_2 P^2$ be a mixture distribution on $\mathbb{R}^l$
with $\pi_1 + \pi_2 = 1$. Given a positive semi-definite kernel
$K$, denote by $\lambda_0$, $\lambda_0^1$ and $\lambda_0^2$ and $\phi_0$, $\phi_0^1$ and $\phi_0^2$, the top eigenvalues and the associated eigenvectors of $K_p$, $K_{p^1}$ and $K_{p^2}$,  respectively. Let $t = \lambda_0-\lambda_1$, the eigenvalue-gap of $K_P$. If the constant $r$ defined above satisfies $r < t$, and
\[
\left\| \pi_2 \int_{\mathbb{R}^d } K(x,y)\phi_0^1(y)dP^2(y) \right\|_{L_P^2} \leq \epsilon
\]
such that $\epsilon  + r < t $, then $\pi_1 \lambda_0^1$ is close to $K_P$'s top eigenvalue $\lambda_0$,
\[
|\pi_1 \lambda_0^1  - \lambda_0 | \leq \epsilon
\]
and $\phi_0^1$ is close to $K_P$'s top eigenfunction $\phi_0$ in the $L_P^2$ sense,
\[ 
\| \phi_0^1 - \phi_0 \|_{L_P^2}	\leq  \frac{\epsilon}{ t-\epsilon}.
\]
\end{corollary}

%\begin{theorem}[Theorem 1, Shi 2009]
%  Given a $d$-dimensional mixture of two Gaussian $p(x) = \sum_{i=1}^2 \pi_i p_i(x)$ where $\pi_i$ is mixing weight and $p_i(x)$ is the density corresponding to $N(\mu_i,\sigma^2 I)$. Define $\beta = 2\sigma^2/\omega^2$ and $\xi = \sqrt{2}\sigma / \sqrt{\sqrt{1+2\beta}-1}$, then the first eigenfunction ($\phi_0^1$ with eigenvalue $\lambda_0^1$) of $G_{p_1}^\omega$ is approximately an eigenfunction of $G_P^\omega$ in the following sense.

%For any $\epsilon > 0$ we have that for all $y$
%\[
%G_P^\omega \phi_0^1 (y)	= \pi_1 \lambda_0^1 [\phi_0^1 (y) + T(y)]	\quad \textrm{and}\quad	|T(y)|\leq \epsilon,
%\]
%assuming that the separation satisfies
%\begin{align*}
%\frac{ \| \mu_1-\mu_2 \|^2 }{\xi^2 + \sigma^2}	\geq 2\log\left( \frac{\pi_1}{\pi_2} \right)	+2 \log\left( \frac{1}{\epsilon} \right)	+\frac{d}{4} \log(1+2\beta).
%\end{align*}
%\end{theorem}

\subsection{Gaussian components}

The eigenvalues and the eigenvectors of the kernel operators for the univariate Gaussian  $N(\mu,\sigma^2)$ were calculated in \cite{Zhu} and \cite{Shi1}. Specifically, if one considers the kernel
\[
	K(x,y) = e^{-\frac{(x-y)^2}{2\omega^2}},\quad
	\omega\in\mathbb{R}_+,
\]
and the corresponding integral operator
\begin{align*}
K_P^\omega f(x) & = \int_{\mathbb{R}} K(x,z) f(z)p(z) dz\\
				& = \int_{\mathbb{R}}	e^{ - \frac{(z-x)^2}{2\omega^2} }f(z)p(z)dz,
\end{align*}
then the following holds.

\begin{theorem}[\cite{Shi1}, Proposition 1] \label{alytic}
Let $\beta = 2\sigma^2/\omega^2$ and let $H_i(x)$ be the $i$-th order Hermite polynomial. Then the eigenvalues and eigenfunctions of $K_P^\omega$ for $i = 0,1,\ldots,n$ are given by
\begin{gather*}
	\lambda_i	= \frac{\sqrt{2}}{(1+\beta+\sqrt{1+2\beta})^{1/2}} 	\left( \frac{\beta}{1+\beta+\sqrt{1+2\beta}}\right)^i,\\
	\phi_i(x) 	= \frac{ (1+2\beta)^{1/8} }{\sqrt{2^i i!}} 	e^{ -\frac{ (x-\mu)^2 }{2\sigma^2} \frac{\sqrt{1+2\beta} 	-1}{2} } H_i \left[ \left( \frac{1+2\beta}{4}\right)^{1/4} \frac{x-\mu}{\sigma} \right].
\end{gather*}
\end{theorem}
The eigenvalues form a geometric series with a common ratio
$${ \beta \over 1+\beta+\sqrt{1+2\beta} } < 1.$$ 
It is clear that the sequence of the eigenvalues converges to zero faster for smaller values of $\beta$.

\subsection{Multivariate Gaussian distribution}

\begin{theorem}[\cite{Shi1}]
  Let $N(\mu,\Sigma)$ be a multivariate Gaussian in $\mathbb{R}^d$. Let $\Sigma = \Sigma_{i=1}^d \sigma_i^2 u_i u_i^t$ be the spectral decomposition of the covariance matrix $\Sigma$. Set
  $$K_P^\omega f(x) = \int_{\mathbb{R}^d} e^{-{\|x-z\|^2 \over 2\omega^2}} f(z)p(z) dz.$$
  Then $K_P^\omega$ can be decompose as a direct sum
  \[
  K_P^\omega = \oplus_{i=1}^d K_{P^i}^\omega
  \]
 where $P_i$ is the one-dimensional Gaussian distribution with variance $\sigma_i^2$ and mean $\langle \mu, u_i \rangle$ along the direction of $u_i$.

Then the eigenvalues and eigenfunctions of $K_P^\omega $ can be written as
\begin{align*}
  \lambda_{[i_1,\ldots,i_d]} &= \prod_{j=1}^d \lambda_{i_j} (K_{p_j}^\omega),\\
  \phi_{[i_1,\ldots,i_d]} (K_P^\omega)(x)&= \prod_{j=1}^d (K_{p_j}^\omega)(\langle x,u_j\rangle) 
\end{align*}
where $[i_i,\ldots,i_d]$ is a multi-index over all components.
\end{theorem}

\section{Estimates on the second eigenvalue of a Gaussian mixture} \label{section_gap}

In this Section we will provide bounds on the second eigenvalue of the kernel operator for a Gaussian mixture $P=\pi_1 P^1+ \pi_2 P^2$. Our goal is to come up with computable bounds which can guarantee, if the parameters of the Gaussian mixture are chosen appropriately, that there is definitely a spectral gap between the top eigenvalue (see Theorem $\ref{shi}$) and the second one. As we will demonstrate with numerical simulation, clustering of a Gaussian mixture with such a spectral gap invariably results in a formation of a single cluster.

\begin{mthm}\label{second_ev}
Let $P=\pi_1 +P^1 + \pi_2 P^2$ be a mixture distribution. Let the top and the second eigenvalue of  $K_P$, $K_{P_1}$ and $K_{P_2}$ 
be denoted by $\lambda_0$, $\lambda_1$,  $\gamma_0$, $\gamma_1$ and $\nu_0$, $\nu_1$, respectively. 
Denote $\|\cdot\|_P = \|\cdot\|_{L^2_P}$, and similarly for $P^1$ and $P^2$. Assume $\pi_1 > \pi_2$, and let
\[
 \delta(z) = \phi_0^1(z)-\phi_0(z).
\]
Then an upper bound for $\lambda_1$ is given by
\begin{align*}
\lambda_1  \leq  \pi_1^2 \gamma_1 \left( \frac{1}{\pi_1} + \frac{2}{\sqrt{\pi_1}} A + A \right)& + 2\pi_1^2  \left\|  K \right\|_{{P^1} \times {P^1}}  \left( \frac{1}{\sqrt{\pi_1}} + A 
\right)A\\
& + \pi_1^2 \left\|  K \right\|_{{P^1} \times {P^1}}  A^2 + r
\end{align*}
where
\begin{align*}
  A \leq \| \delta\|^2_{P}  + 2 \| \delta\|_P \Delta 
  + \Delta^2, \quad \Delta &\leq \left(\frac{1}{\pi_1}-1\right) \| \phi_0^1 \|_P  +
\frac{1}{\pi_1} \| \phi_0^1 \|_{P^2}
\end{align*}
and
\[
	r = \left(\pi_1 \pi_2  \iint K(x,z)^2 dP^1(x)dP^2(z)\right)^{1 \over 2}.
\]

Moreover, a lower bound is given by
\begin{equation*}
\begin{split}
	\lambda_1 & \geq \frac{1}{D_1+D_2}\left[ \hspace{-2mm} \phantom{\int} \hspace{-2mm}\pi_1^2 \gamma_1+\pi_1 \pi_2 \|\phi_1^1\|^2_{P^2}-\pi_2 \|\phi_1^1\|_P\|\phi_1^1\|_{P^2} \|K\|_{P^2\times P}- \right.\\
	&\hspace{30pt} \left. -2|e|\left|
	 \lambda_0 \cdot
		\| \delta \|_P\cdot \| \phi_1^1\|_P + \lambda_0 \int \phi_1^1(x)\phi_0^1(x) dP(x) \right|-e^2\lambda_0\right]
\end{split}
\end{equation*}
where
\begin{align*}
	D_1 & = \pi_1 + 2|e|\pi_1 \|\delta\|_{P} + e^2 + 2e^2\pi_1  \|\delta \|_{P} + e^2 \pi_1  \|\delta\|^2_{P},\\
	D_2 & = \pi_2 \|\phi_1^1\|^2_{P^2} + 2|e| \pi_2  \|\delta\|_{P} \|\phi_1^1\|_{P^2} + e^2 \pi_2\|\phi_0^1\|^2+ \\
        \nonumber & \phantom{ = \pi_2 \|\phi_1^1\|^2_{P^2} } + 2e^2\pi_2  \|\delta \|_{P}\|\phi_0^1\|_{P^2} + e^2 \pi_2  \|\delta\|^2_{P}
\end{align*}
and where
\[
	e = \int \phi_1^1(\phi_0^1-\phi_0)dP - \pi_2 \int \phi_1^1 \phi_0^1 dP^2.
\]
\label{the_theorem}
\end{mthm}

\begin{remark}
If we assume Corollary \ref{cor}, then
\[
\|\delta \|_P \leq \epsilon/(t-\epsilon).
\]
\end{remark}

\begin{remark}
All norms in the above theorem are estimated in the Appendix.
\end{remark}
\begin{proof} \textbf{Upper bound.} We first prove the upper bound. For this, we present the standard calculation for $f$ in the othogonal compliment $V$ of the span of $\phi_0$ in $L^2_P$. We have
\begin{equation}\label{lambda1}
  \lambda_1 = \max_{f \in  V}  \frac{\iint K(x,z) f(x)f(z) dP(x)dP(z) }{\int f(x)^2 dP(x)},
\end{equation}
and that for any such  $f$,

\begin{align}
 \iint & K(x,z) f(x)f(z) dP(x)dP(z)  = \pi_1^2 \iint K(x,z)f(x)f(z)dP^1(x) dP^1(z) \notag \\
 & \hspace{20pt}+ \pi_2^2 \iint K(x,z)f(x)f(z)dP^2(x)dP^2(z)+ \notag \\
 &\hspace{20pt}+2\pi_1 \pi_2 \iint K(x,z)f(x)f(z)dP^1(x)dP^2(z) \notag \\
 & \le  \pi_1^2  \iint K(x,z) \left(f(z) - (\textrm{proj}^1_{\phi_0^1} f)\phi_0^1(z) +F(z)\right)  \times \notag \\
 &\hspace{20pt} \times \left(f(x) - (\textrm{proj}^1_{\phi_0^1} f)\phi_0^1(x) +F(x)\right)  dP^1(z)  d P^1(x) + \notag \\
&\hspace{20pt}+ \pi_2^2 \nu_0 \int (f(x))^2 dP^2(x) +\notag \\
 &\hspace{20pt}  + 2\pi_1 \pi_2 \iint K(x,z)f(x)f(z)dP^1(x)dP^2(z), \label{upper1} 
\end{align}
where
$$\textrm{proj}^i_{\phi} f = \int \phi(x)f(x) dP^i(x), \quad  \textrm{proj}_{\phi} f = \int \phi(x)f(x) dP(x)$$
and
$$F(z)  = (\textrm{proj}^1_{\phi_0^1} f)\phi_0^1(z)= (\textrm{proj}^1_{\phi_0^1} f - \textrm{proj}_{\phi_0} f)\phi_0^1(z)$$
and where the second equality follows from $f\in V$. Therefore, the bound $(\ref{upper1})$ becomes
\begin{eqnarray} \nonumber   \iint \hspace{-3mm}  &\hspace{-3mm} & \hspace{-3mm} K(x,z) f(x)f(z) dP(x)dP(z)  \\
 \nonumber &&   \le  \pi_1^2 \gamma_1 \int  \left(f(x) - (\textrm{proj}^1_{\phi_0^1} f)\phi_0^1(x) \right)^2   d P^1(x) +\\
\nonumber &&  \hspace{10pt}+ 2 \pi_1^2 \iint  K(x,z) \left(f(z) - (\textrm{proj}^1_{\phi_0^1} f)\phi_0^1(z) \right)  F(x)  \ dP^1(z)  d P^1(x) +\\
\nonumber &&  \hspace{10pt}+ \pi_1^2  \iint  K(x,z)  F(z) F(x)  \ dP^1(z)  d P^1(x) +\\
\nonumber && \hspace{10pt} + \pi_2^2 \nu_0 \int f(x)^2 dP^2(x)+ \\
\label{upper2} && \hspace{10pt} +2\pi_1 \pi_2 \iint K(x,z)f(x)f(z)dP^1(x)dP^2(z). 
\end{eqnarray}
We take the terms one by one in (\ref{upper2}). For the first, we write
\begin{align*}
\int \left( f(x)-(\textrm{proj}^1_{\phi_0^1} f)\phi_0^1(x) \right)^2 dP^1(x)
& = \|f-(\textrm{proj}^1_{\phi_0^1})\phi_0^1\|^2 \\
& \leq \left( \|f\|_{P^1} + \|F\|_{P^1} \right)^2  \\
& \leq \left( \frac{1}{\sqrt{\pi_1}} + \|F\|_{P^1} \right)^2
\end{align*}
where we in the last step have normalized $f$. 

For the second term, we write
\begin{align*}
& \left|\iint K(x,z)\left(f(z)-(\textrm{proj}_{\phi_0^1} f) \phi_0^1(z)\right)F(x)dP^1(z)dP^1(x) \right|\\
& \hspace{20pt} =  \left| \int \left( \int K(x,z)\left(f(z)-(\textrm{proj}_{\phi_0^1} f)
\phi_0^1(z) \right) dP^1(z) \right) F(x)  dP^1(x) \right|\\
& \hspace{20pt} \le  \sqrt{ \int \left( \int K(x,z)\left( f(z)-(\textrm{proj}_{\phi_0^1} f)\phi_0^1(z) \right) dP^1(z) \right)^2 dP^1(x) } \left\| F  \right\|_{{P^1}}  \\
  & \hspace{20pt} \le  \sqrt{ \int  \int K(x,z)^2 d P^1(z) \left\|f-(\textrm{proj}_{\phi_0^1} f)\phi_0^1  \right\|^2_{{P^1}}  dP^1(x) } \left\| F  \right\|_{{P^1}}  \\
  & \hspace{20pt} \le \left\|  K \right\|_{{P^1} \times {P^1}}  \left\|f-(\textrm{proj}_{\phi_0^1} f)\phi_0^1  \right\|_{{P^1}} \left\| F  \right\|_{{P^1}}\\
   & \hspace{20pt} \le \left\|  K \right\|_{{P^1} \times {P^1}}  
   \left( \| f\|_{P^1} + \|F\|_{P^1} \right)
   \left\| F  \right\|_{{P^1}}\\
   & \hspace{20pt} \le \left\|  K \right\|_{{P^1} \times {P^1}}  
      \left( \frac{1}{\sqrt{\pi_1}} + \|F\|_{P^1} \right)
      \left\| F  \right\|_{{P^1}}.
\end{align*}
%We refer to the Appendix for calculating 
%$\left\|  K \right\|_{{P^1} \times {P^1}}$. 

For the third term,
$$ \left|\iint K(x,z) F(z) F(x) dP^1(z)dP^1(x) \right| \le \left\|  K \right\|_{{P^1} \times {P^1}}  \left\| F  \right\|^2_{{P^1}}.$$
%We refer to the Appendix for calculating $\|K\|_{P^1\times P^1}$.

For the fourth,
\begin{align*}
\pi_2^2 \nu_0 \int f(x)^2 dP^2(x) \leq \pi_2 \nu_0.
\end{align*}

For the last term,
\begin{align*}
  & 2\pi_1 \pi_2 \iint K(x,z) f(x)f(z) dP^1(x)dP^2(z)\\
  & \hspace{5pt} \leq 2\pi_1 \pi_2 \sqrt{ \iint K(x,z)^2 dP^1(x)dP^2(z) } \sqrt{ \iint f(x)^2f(z)^2 dP^1(x)dP^2(z) }\\
  & \hspace{5pt} = 2\sqrt{\pi_1\pi_2 \iint K(x,z)^2  dP^1(x) dP^2(z) } \times \\
  & \hspace{140pt} \times \sqrt{ \pi_1 \int f(x)^2 dP^1(x) } \sqrt{ \pi_2 \int f(z)^2 dP^2(x) }\\
  & \hspace{5pt} \leq \sqrt{\pi_1\pi_2 \iint K(x,z)^2  dP^1(x) dP^2(z) } \times \\
  & \hspace{140pt} \times \left(  \pi_1 \int f(x)^2 dP^1(x) +  \pi_2 \int f(z)^2 dP^2(x) \right)\\
& \hspace{5pt} = r\int f(x)^2 dP(x),
\end{align*}
where
$$r = \left(\pi_1 \pi_2  \iint K(x,z)^2 dP^1(x)dP^2(z)\right)^{1 \over 2}.$$

Finally,
\begin{align*}
\lambda_1 & \leq  \pi_1^2 \gamma_1 \left( \frac{1}{\sqrt{\pi_1}} + \| F \|_{P^1} \right)^2+ 2\pi_1^2 \left( \left\|  K \right\|_{{P^1} \times {P^1}}  \left( \frac{1}{\sqrt{\pi_1}} + \|F\|_{P^1}  \right)\left\| F  \right\|_{{P^1}}\right)\\
& \phantom{  \leq  \pi_1^2 \gamma_1 \left( \frac{1}{\sqrt{\pi_1}} +   \| F \|_{P^1} \right)^2 \ } +   \pi_1^2 \left\|  K \right\|_{{P^1} \times {P^1}}  \left\| F  \right\|^2_{{P^1}} + r.
\end{align*}

\noindent \textbf{Lower bound.} Let
$$\widetilde{\phi}_1^1 =\phi_1^1 - (\textrm{proj}_{\phi_0}\phi_1^1)\phi_0= \phi_1^1 + \phi_0(\textrm{proj}^1_{\phi_0^1} \phi_1^1-\textrm{proj}_{\phi_0}\phi_1^1),$$
where in the second equality we have used that $(\textrm{proj}^1_{\phi_0^1} \phi_1^1)\phi_0 = 0$.
Further,
\begin{align*}
\textrm{proj}^1_{\phi_0^1} \phi_1^1- \textrm{proj}_{\phi_0}\phi_1^1
& = \int \phi_0^1(z)\phi_1^1(z) dP^1(z) - \int \phi_0(z)\phi_1^1(z) dP(z)\\
& = \pi_1 \int \phi_0^1(z)\phi_1^1(z) dP^1(z) - \int \phi_0(z)\phi_1^1(z) dP(z)\\
& =  \int \phi_0^1(z)\phi_1^1(z) dP(z) - \int \phi_0(z)\phi_1^1(z) dP(z)\\
& \hspace{30pt} - \pi_2 \int \phi_1^1(z)\phi_0^1(z) dP^2(z)\\
& =  \int \phi_1^1(z)(\phi_0^1(z) - \phi_0(z)) dP(z)\\
& \hspace{30pt} - \pi_2 \int \phi_1^1(z)\phi_0^1(z) dP^2(z).
\end{align*}

Let now
\begin{align}
  E(x)& = \phi_0(x) \left(  \int \phi_1^1(z) (\phi_0^1(z)-\phi_0(z))  dP(z) - \pi_2 \int \phi_1^1(z) \phi_0^1(z) dP^2(z)\right) \notag \\
  &=e \phi_0(x), \notag 
\end{align}
where
$$e= \int \phi_1^1(z) (\phi_0^1(z)-\phi_0(z))  dP(z) - \pi_2 \int \phi_1^1(z) \phi_0^1(z) dP^2(z),$$
so that, by definition,
\begin{align}
\nonumber \lambda_1 & \geq  \frac{\iint K(x,z)\widetilde{\phi}_1^1(x)\widetilde{\phi}_1^1(z) dP(x) dP(z) }{ \int (\widetilde{\phi}_1^1(x))^2 dP(x) } \\
\nonumber & =  \frac{\iint K(x,z) (\phi_1^1(x)+E(x))(\phi_1^1(z)+E(z)) dP(x) dP(z) }{ \int [\phi_1^1(x)+E(x)]^2 dP(x) }\\
\nonumber & =  \frac{\iint K(x,z)\phi_1^1(x)\phi_1^1(z) dP(x) dP(z) }{ \int [\phi_1^1(x)+E(x)]^2 dP(x) }+ \frac{\iint K(x,z)\phi_1^1(x)E(z) dP(x)
dP(z) }{ \int [\phi_1^1(x)+E(x)]^2 dP(x) } \\
& \phantom{  =  \frac{\iint K(x,z)\phi_1^1(x)\phi_1^1(z) dP(x) dP(z) }{ \int [\phi_1^1(x)+E(x)]^2 dP(x) } \ }+   \frac{\iint K(x,z)E(x)\phi_1^1(z) dP(x) dP(z) }{ \int [\phi_1^1(x)+E(x)]^2 dP(x) } \nonumber \\
& \phantom{  =  \frac{\iint K(x,z)\phi_1^1(x)\phi_1^1(z) dP(x) dP(z) }{ \int [\phi_1^1(x)+E(x)]^2 dP(x) } \ } + \frac{\iint K(x,z)E(x)E(z) dP(x) dP(z) }{ \int [\phi_1^1(x)+E(x)]^2 dP(x) }.
\label{lower}
\end{align}
We will  bound terms in the {\bf numerator} of $(\ref{lower})$  one by one.

\medskip

\noindent {\it First term in $(\ref{lower})$}. The numerator in the first term is
\begin{align}
&\iint K(x,z)\phi_1^1(x)\phi_1^1(z)dP(x)dP(z)\nonumber\\
&\hspace{20pt} = \int \left(\int K(x,z)\phi_1^1(x)dP(x)\right)\phi_1^1(z) dP(z)\nonumber\\
&\hspace{20pt} = \int \left( \int \pi_1 K(x,z)\phi_1^1(x)dP^1(x)
+ \int \pi_2 K(x,z)\phi_1^1(x)dP^2(x)\right)\phi_1^1(z) dP(z)\nonumber\\
&\hspace{20pt} = 
\int\left( \pi_1 \gamma_1 \phi_1^1(z) + 
\int \pi_2 K(x,z)\phi_1^1(x)dP^2(x)\right)\phi_1^1(z) dP(z)\nonumber\\
&\hspace{20pt} = 
\int \pi_1 \gamma_1 \phi_1^1(z)\phi_1^1(z) dP(z)
+ \iint \pi_2 K(x,z)\phi_1^1(x)\phi_1^1(z) dP^2(x)dP(z).\nonumber
\end{align}

We take the terms one by one, and write 
\begin{align*}
& \int \pi_1 \gamma_1\phi_1^1(z)\phi_1^1(z) dP(z)\\
& \hspace{20pt} = \int  \pi_1^2 \lambda_1\phi_1^1(z)\phi_1^1(z)dP^1(z)
+
\int \pi_1 \pi_2  \gamma_1\phi_1^1(z)\phi_1^1(z)dP^2(z)\\
& \hspace{20pt} = \pi_1^2 \gamma_1 + \gamma_1\pi_1\pi_2 \| \phi_1^1 \|^2_{P^2}.
\end{align*}

For the second, we write
\begin{eqnarray}\label{omskr}
\nonumber && \iint \pi_2 K(x,z)\phi_1^1(x)\phi_1^1(z) dP^2(x)dP(z)  \\ 
\nonumber && \hspace{10pt} = \int \phi_1^1(z)\left( \int \pi_2 K(x,z)\phi_1^1(x)dP^2(x)\right) dP(z) \\
\nonumber && \hspace{10pt} \leq  \sqrt{\int \phi_1^1(z)^2 dP(z) }  \sqrt{ \int \left(\int \pi_2 K(x,z)\phi_1^1(x)dP^2(x)\right)^2 dP(z)} \\
\nonumber && \hspace{10pt} \leq \pi_2  \sqrt{\int \phi_1^1(z)^2 dP(z) } \sqrt{ \int  \int K(x,z)^2 dP^2(x)  \int \phi_1^1(x)^2 dP^2(x) \  d P(z)} \\
 && \nonumber \hspace{10pt} \leq \pi_2  \left\| \phi^1_1 \right\|_{P}  \left\| \phi_1^1 \right\|_{P^2}  \left\| K \right\|_{P^2 \times {P}}.
\end{eqnarray}
%We refer to the Appendix for calculating $\left\| \phi_1^1 \right\|_{{P^2}}$ and $\left\| K \right\|_{{P^2} \times {P}}$.\\

\noindent {\it Second term in $(\ref{lower})$}. We write
\begin{equation*}
%\begin{split}
 \iint K(x,z)\phi_1^1(x)E(z)dP(x)dP(z)= e \!\! \iint K(x,z)\phi_0(z)\phi_1^1(x) dP(x)dP(z)
%\end{split}
\end{equation*}
%where $e = \int \phi_1^1(\phi_0^1-\phi_0)dP - \pi_2 \int \phi_1^1 \phi_0^1 dP^2$.
We focus on
\begin{align}
	& \iint K(x,z) \phi_0(z)\phi_1^1(x)dP(x)dP(z)\nonumber\\
	&\hspace{20pt} = \int \phi_1^1(x) \left(
	\int K(x,z)\phi_0(z) dP(z)\right) dP(x)\nonumber\\
	&\hspace{20pt} =
	\int \phi_1^1(x) \lambda_0 \phi_0(x) dP(x)\nonumber\\
	&\hspace{20pt}
	= \lambda_0 \int \phi_1^1(x)(\phi_0^1(x)+\delta(x)) dP(x) \nonumber\\
	&\hspace{20pt} = \lambda_0 \int \delta (x)\phi_1^1(x)dP(x)
	+\lambda_0 \int \phi_1^1(x)\phi_0^1(x) dP(x)\nonumber\\
	&\hspace{20pt} \leq \lambda_0 \cdot
	\| \delta \|_P\cdot \| \phi_1^1\|_P + \lambda_0 \int \phi_1^1(x)\phi_0^1(x) dP(x).
\label{h}
\end{align} 
%We refer to the Appendix for calculating $\int \phi_1^1 \phi_0^1 dP$.\\

Since the {\it Third term in $(\ref{lower})$} is precisely symmetrical to the second, we go directly to the {\it Fourth term in $(\ref{lower})$}. Here we have
\begin{align*}
&e^2 \iint K(x,z)\phi_0(x)\phi_0(z) dP(x)dP(z) = e^2 \lambda_0.
\end{align*}

For the {\bf denominator of $(\ref{lower})$}  we write
\begin{align*}
	\int (\phi_1^1(x) + E(x))^2 dP(x) & = \int (\phi_1^1(x)+\phi_0(x)e)^2dP(x)\\
	& = \int (\phi_1^1(x)+[\phi_0^1(x)+\delta(x)]e)^2 dP(x).
\end{align*}

We focus first on integrating with respect to $P^1$,
\begin{align}
  \nonumber & \pi_1 \int (\phi_1^1(x)+[\phi_0^1(x)+\delta(x)]e)^2dP^1(x)
  \\
  \nonumber &\hspace{30pt} = \pi_1\int \phi_1^1(x)^2 dP^1(x) + 2\pi_1 \int \phi_1^1(x)[\phi_0^1(x)+\delta(x)]e dP^1(x)
  \\
  \nonumber &\hspace{50pt} + e^2\pi_1 \int (\phi_0^1(x)+\delta(x))^2 dP^1(x)
  \\
  \nonumber &\hspace{30pt} = \pi_1 + 2 e \pi_1  \int \phi_1^1(x)\delta(x) dP^1(x)+ e^2 \pi_1 \int (\phi_0^1(x))^2 dP^1(x) \\
 \nonumber &\hspace{50pt}  + 2e^2\pi_1  \int \phi_0^1(x)\delta(x)dP^1(x)
  + e^2\pi_1  \int (\delta(x))^2 dP^1(x) \\
  \nonumber &\hspace{30pt} \leq \pi_1 + 2|e|\pi_1  \|\delta\|_{P^1} + e^2 \pi_1  + 2e^2 \pi_1  \|\delta \|_{P^1} + e^2 \pi_1  \|\delta\|^2_{P^1} \\
  &\hspace{30pt} \leq \pi_1 + 2|e|\pi_1 \|\delta\|_{P} + e^2\pi_1 + 2e^2\pi_1  \|\delta \|_{P} + e^2 \pi_1  \|\delta\|^2_{P}.  \label{firstterm}
\end{align}
Similarly, for $P^2$,
\begin{align}
  \nonumber	&\hspace{40pt} \pi_2 \int (\phi_1^1(x)+[\phi_0^1(x)+\delta(x)]e)^2dP^2(x) \\
  \nonumber	&\hspace{70pt} = \pi_2\int (\phi_1^1(x))^2 dP^2(x) 
  + 2\pi_2 \int \phi_1^1(x)[\phi_0^1(x)+\delta(x)]e dP^2(x) \\
  \nonumber	&\hspace{90pt} + e^2\pi_2 \int (\phi_0^1(x)+\delta(x))^2 dP^2(x) \\
  \nonumber	&\hspace{70pt} = \pi_2\int (\phi_1^1(x))^2 dP^2(x) + 2e\pi_2  \int \phi_1^1(x)\delta(x) dP^2(x) \\
  \nonumber	&\hspace{90pt} + e^2 \pi_2 \int (\phi_0^1(x))^2 dP^2(x)+ 2e^2\pi_2  \int \phi_0^1(x)\delta(x)dP^2(x) \\
  \nonumber     &\hspace{90pt}	+ e^2\pi_2  \int (\delta(x))^2 dP^2(x)\\
  \nonumber	&\hspace{70pt} \leq \pi_2 \|\phi_1^1\|^2_{P^2} + 2|e| \pi_2  \|\delta\|_{P} \|\phi_1^1\|_{P^2} + e^2 \pi_2\|\phi_0^1\|^2_{P^2} \\
  &\hspace{90pt} + 2e^2\pi_2  \|\delta \|_{P}\|\phi_0^1\|_{P^2} + e^2 \pi_2  \|\delta\|^2_{P}. \label{secondterm}
\end{align}
%We refer to the Appendix in calculating $\|\phi_0^1\|^2_{P^2}$.
\end{proof}

\section{Clustering and simulations}\label{sims1}

\subsection{Single component in $\mathbb{R}$}

\begin{figure}[h!]%
    \centering
\begin{tabular}{c c}
  {{\includegraphics[width=0.49\textwidth]{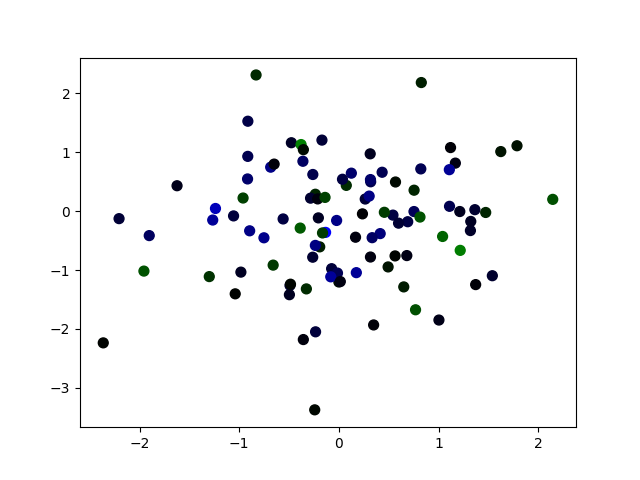} }}%
  & {{\includegraphics[width=0.49\textwidth]{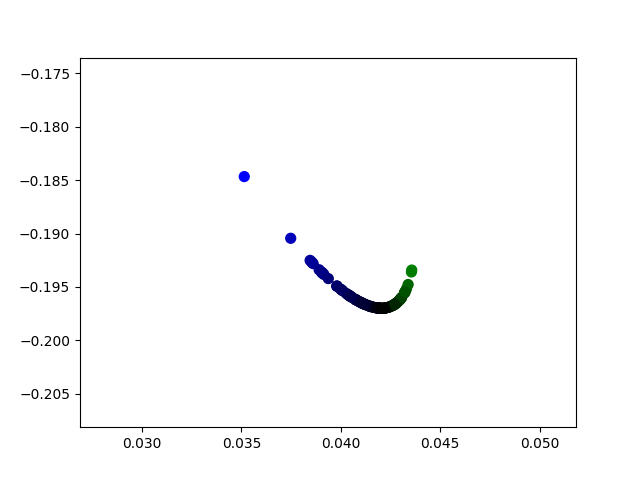} }}%
\end{tabular}
    \caption{Clustering in $\mathbb{R}$. Single component. $\mu = 10$,
    $\sigma^2 = 1$. 
    The algorithm clusters the points in a continuous spectrum measured
    with respect to color.}%
    \label{old_small}%
\end{figure}

Fig. \ref{old_small} presents clustering
after one iteration. The parameters are set to
$\mu = 10$, $\sigma^2 = 1$.
After repeated iterations
the points tend either to contract to a single
point or form a straight line. In general
the clustering seems satisfactory given that the deviation
is not too high. 
In a certain sense, the clustering \emph{respects
the underlying data structure}: if there, to begin with, 
do exist 
clusters in the underlying data, these
will show up in the final plot. And conversely: if there are no clusters
to begin with (the points may be very dispersed because of a high deviation), no clusters will later show up.

\subsection{Mixture of two components in $\mathbb{R}$}

\begin{figure}[h!]%
  \centering
  \begin{tabular}{c c}
    {{\includegraphics[width=0.49\textwidth]{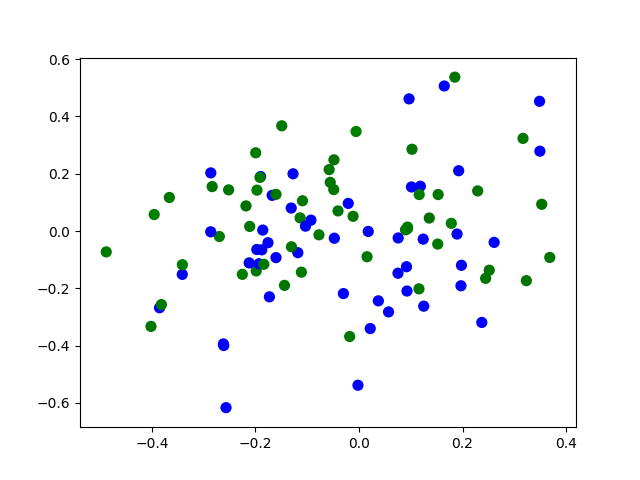} }}%
    & {{\includegraphics[width=0.49\textwidth]{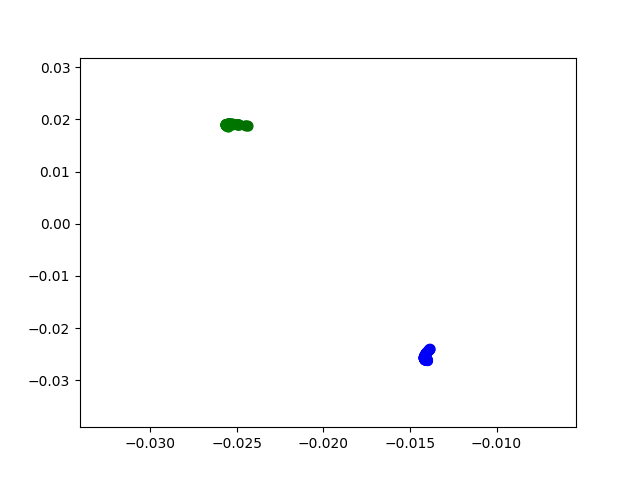} }}%
  \end{tabular}
  \caption{Clustering in $\mathbb{R}$. Two components.
    $\mu_1 = -5$, $\mu_2 = 5$, $\sigma_1^2 = 0.1$, $\sigma_2^2 = 0.1$.
    Because of the high separation of distributions, two distinct and well-separated clusters are produced. The algorithm  preserves, or, respects the underlying data structure.}%
  \label{mix_two}%
\end{figure}

Fig. \ref{mix_two} presents a clustering of a mixture with two components.
The parameters are set to
$\mu_1 = -5$, $\mu_2 = 5$, $\sigma_1^2 = 0.1$, $\sigma_2^2 = 0.1$.
Note the high separation
of distributions, accurately portrayed in the final clustering by the large
separation between the two blue and green groups (the distance \emph{between}
groups is large compared to the distance \emph{within} each group).

\subsection{Single component in $\mathbb{R}^2$}
\begin{figure}[h!]%
  \centering
  \begin{tabular}{c c}
    {{\includegraphics[width=0.49\textwidth]{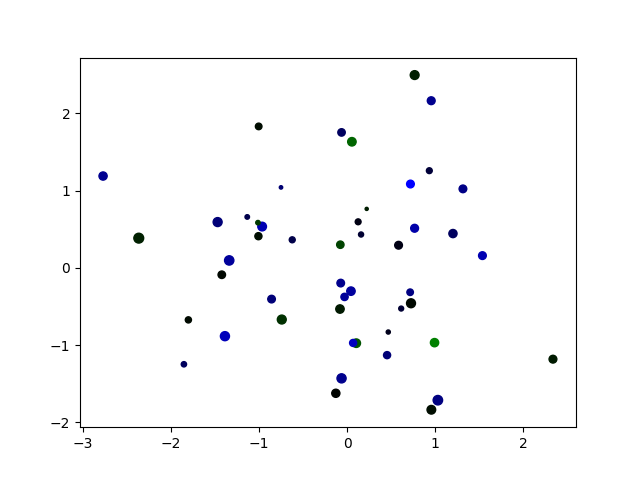} }}%
    & {{\includegraphics[width=0.49\textwidth]{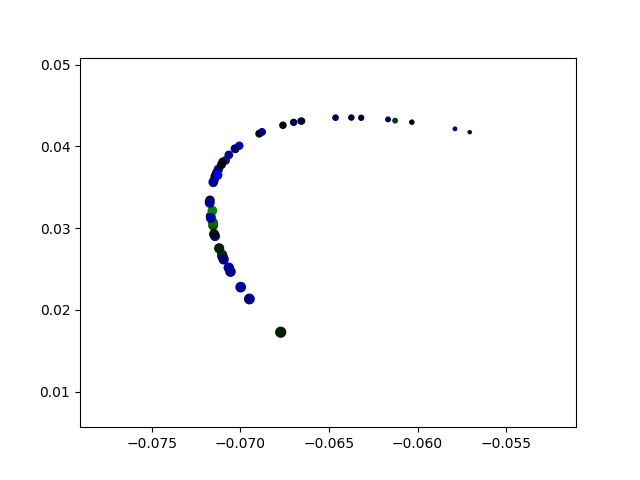} }}%
  \end{tabular}
  \caption{Clustering in $\mathbb{R}^2$. Single component.
    The color parameters are set to $\mu_1 = 20$, $\sigma_1^2 = 1$.
    The size parameters are set to
    $\mu_2 = 10$, $\sigma_2^2 = 0.1$.
    The weights are set to $\alpha_1 = 0.0001$, $\alpha_2 = 1$.
    Because  $\alpha_2 \gg \alpha_1$, the algorithm
    practically discards color in flavor of size.}%
    \label{mul_1}%
\end{figure}

Besides color, one could also measure, for instance, size.
These measures could have weights, reflected in the
kernel
\[
	K_{i, j} = e^{\alpha_1 (c_i-c_j)^2 + \alpha_2 (s_i-s_j)^2}
\]
where $\alpha_1, \alpha_2\in\mathbb{R}_+$.
Fig. \ref{mul_1} presents a clustering using this kernel.
The color parameters are set to $\mu_1 = 20$, $\sigma_1^2 = 1$.
The size parameters are set to $\mu_2 = 10$, $\sigma_2^2 = 0.1$.
Here we have constructed the weights
$\alpha_1$ and $\alpha_2$ so that the algorithm
practically discards color in flavor of size.
In this particular example, $\alpha_1 = 0.0001$ and $\alpha_2 = 1$.

\subsection{Mixture of two components in $\mathbb{R}^2$}

\begin{figure}[h!]%
  \centering
  \begin{tabular}{c c}
    {{\includegraphics[width=0.49\textwidth]{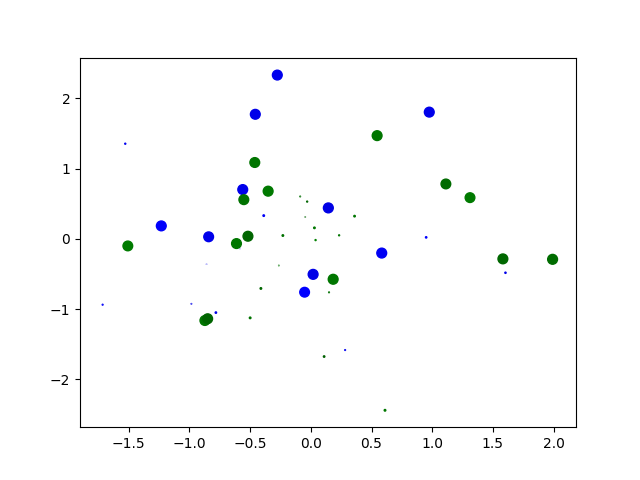} }}%
    &{{\includegraphics[width=0.49\textwidth]{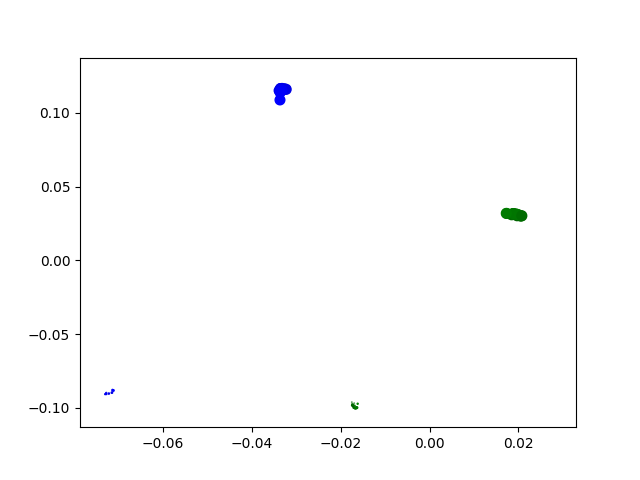} }}%
  \end{tabular}
  \caption{Clustering in $\mathbb{R}^2$. Two components.
    The color parameters are set to 
    $\mu_1 = -20$, $\mu_2 = 20$,
    $\sigma_1^2 = 1$, $\sigma_2^2 = 1$.
    The size parameters are set to
    $\mu_3 = 10$, $\mu_4 = 20$,
    $\sigma_3^2 = 0.1$, $\sigma_4^2 = 0.1$.
    The weights are set to $\alpha_1 = 0.01$ and $\alpha_2 = 1$.
    Four distinct clusters are seen.}%
  \label{mul_3}%
\end{figure}

Fig. \ref{mul_3} presents a mixture 
of two components in $\mathbb{R}^2$.
The color parameters are set to 
$\mu_1 = -20$, $\mu_2 = 20$,
$\sigma_1^2 = 1$, $\sigma_2^2 = 1$.
The size parameters are set to
$\mu_3 = 10$, $\mu_4 = 20$, $\sigma_3^2 = 0.1$, $\sigma_4^2 = 0.1$.
The weights are set to $\alpha_1 = 0.01$, $\alpha_2 = 1$.
Four distinct clusters are seen.

\subsection{Clustering of lines or shapes}

\begin{figure}[h!]%
  \centering
  \begin{tabular}{c c}
    {{\includegraphics[width=0.49\textwidth]{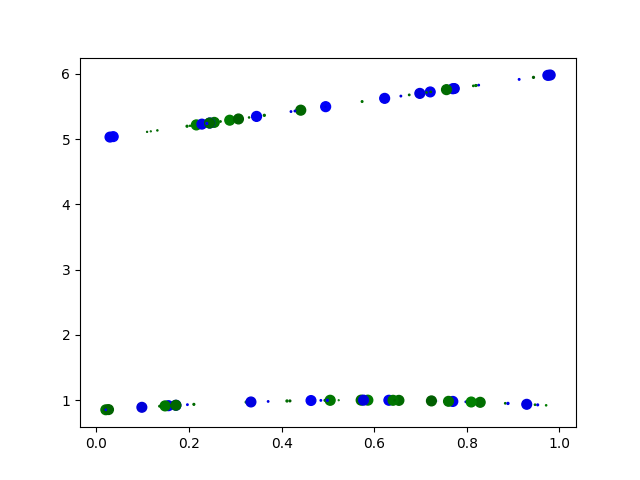} }}%
    &
    {{\includegraphics[width=0.49\textwidth]{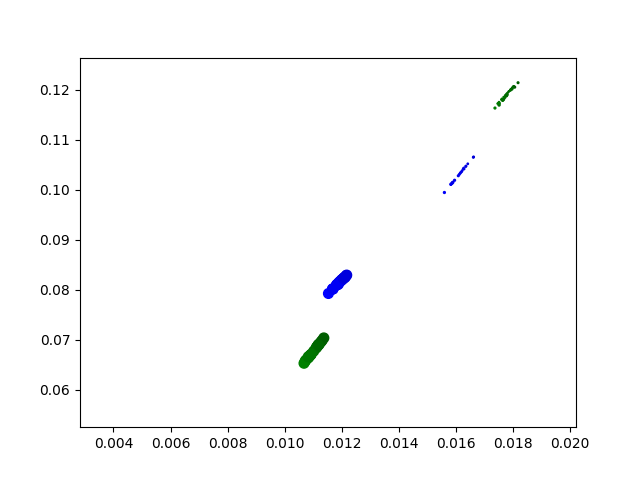} }}%
  \end{tabular}
  \caption{Clustering of lines in $\mathbb{R}^2$. Two components.
    The color parameters are set to
    $\mu_1 = -20$, $\mu_2 = 20$,
    $\sigma_1^2 = 1$, $\sigma_2^2 = 1$.
    The size parameters are set to
    $\mu_3 = 10$, $\mu_4 = 20$,
    $\sigma_3^2 = 0.1$, $\sigma_4^2 = 0.1$.
    The weights are set to $\alpha_1 = 0.01$, $\alpha_2 = 1$.
    Four distinct clusters are seen. Independent of initial structures, the algorithm
    separates the points into distinct and well separated clusters.}%
  \label{lines}%
\end{figure}

The clustering seems (more or less) independent of the initial
shapes. Sometimes objects 
may be placed on lines (underlying clusters themselves), but
the resulting clustering is not affected.
Fig. \ref{lines} presents a mixture 
of two components in $\mathbb{R}^2$, in which the initial points
lie on some given structure.
The color parameters are set to 
$\mu_1 = -20$, $\mu_2 = 20$,
$\sigma_1^2 = 1$, $\sigma_2^2 = 1$.
The size parameters are set to
$\mu_3 = 10$, $\mu_4 = 20$,
$\sigma_3^2 = 0.1$, $\sigma_4^2 = 0.1$.
The weights are set to $\alpha_1 = 0.01$, $\alpha_2 = 1$.
Four distinct clusters are seen.

\subsection{Mixture of two components in $\mathbb{R}$, redefined model}

\begin{figure}[h!]%
  \centering
  \begin{tabular}{c c}
    {{\includegraphics[width=0.49\linewidth]{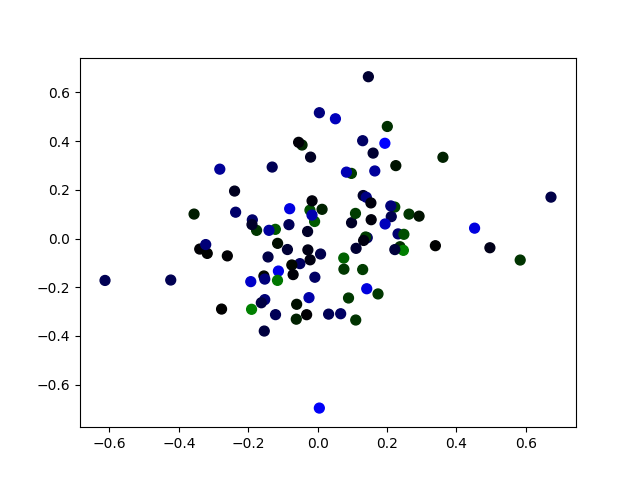} }}%
      &
      {{\includegraphics[width=0.49\linewidth]{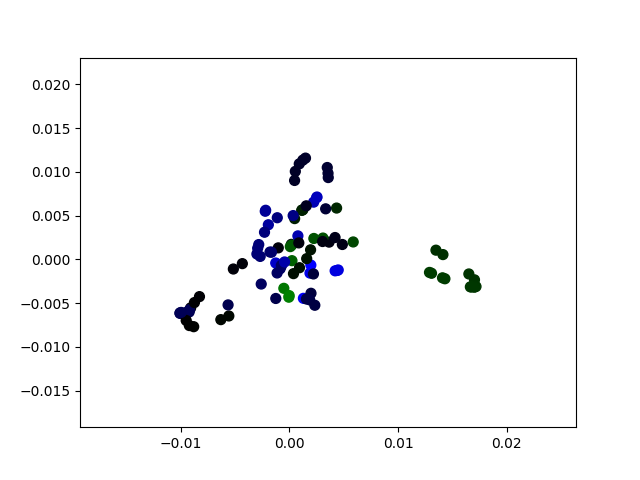} }}%
  \end{tabular}
  \caption{Clustering in $\mathbb{R}$. Two components.
    $\mu_1 = -10$, $\mu_2 = 10$,
    $\sigma_1^2 = 10$, $\sigma_2^2 = 10$.
    Because of the high values of deviations, 
    the points are very much dispersed in feature and will therefore
    not be clustered efficiently by the algorithm.
  }%
  \label{univ_two10}%
\end{figure}

\begin{figure}[h!]%
  \centering
  \begin{tabular}{c c}
    {{\includegraphics[width=0.49\linewidth]{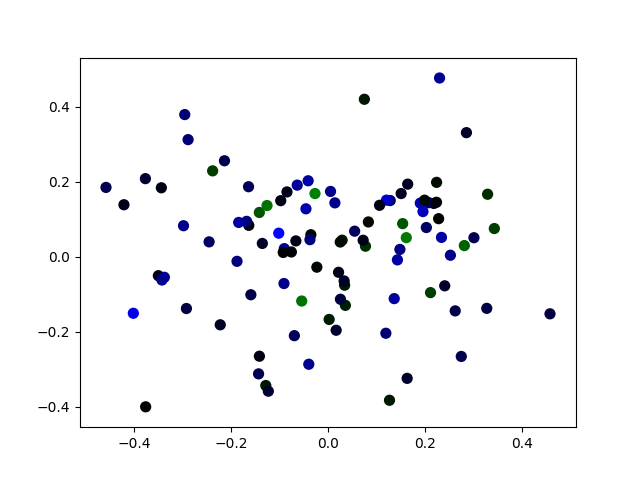} }}%
    &
    {{\includegraphics[width=0.49\linewidth]{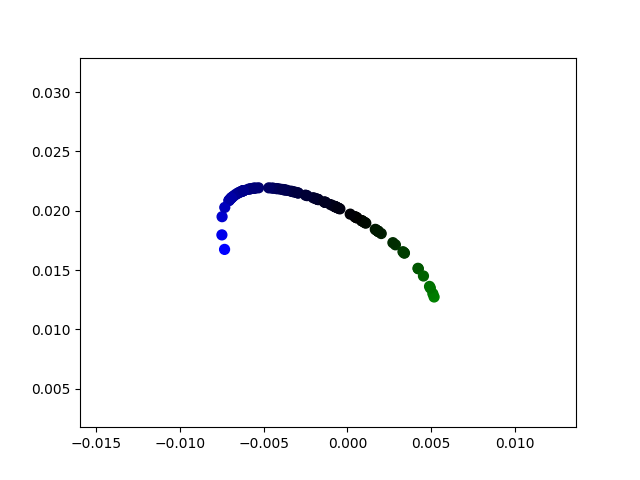} }}%
  \end{tabular}
  \caption{Clustering in $\mathbb{R}$. Two components. Redefined model.
    $\mu_1 = -10$, $\mu_2 = 10$, 
    $\sigma_1^2 = 10$, $\sigma_2^2 = 10$.
    Although being dispersed initially,
    the redefined algorithm produces a remarkably satisfactory clustering.
     }%
  \label{univ_two1}%
\end{figure}

As an adaptation of the original model we suggest the following
redefinition of the kernel,
\begin{align}
\label{red}
	K_{i, j} = \frac{ e^{-\frac{1}{\sigma n} \|x_i-x_j\|^2 / \omega^2} }{n},
\end{align}
where $\sigma = \max(\sigma_1,\ldots,\sigma_l)$.
Clustering with this kernel seems more stable
than the original one.
Fig. \ref{univ_two10} represents a clustering after one iteration
with the original model.
The parameters are set to $\mu_1 = -10$, $\mu_2 = 10$ and
$\sigma_1 = 10$, $\sigma_2 = 10$. Note the high values of deviations, 
making the points very much dispersed in parameter space. 
The effect of scaling by $\sigma n$ is clearly seen if comparing
this last picture with that of Fig. \ref{univ_two1},
where the kernel given in (\ref{red}) was used. 
Although being dispersed initially,
the points seem to cluster remarkably well under this redefined kernel.

%\begin{remark}
%Another suggestion for an improved algorithm is to take into account also the
%\emph{distance between particles} in the kernel. Thus $K_{i, j} = e^{-\|x_i-x_j\|^2-\|y_i-y_j\|^2}$ where $y_i$, $y_j$ are the particle positions. Our experiments show, however, that the result do not improve much as the dynamical system $T_x(y)$ is made non-linear. For this problem we suggest that (\ref{red}) is preferred. 
%\end{remark}

\section{Example of the main theorem}\label{sims2}

In the following section we illustrate the main theorem (see Table \ref{table}
for the relevant parameters). The top three
(experimental) eigenvalues  were calculated to be
$\lambda_0 = 0.62$, $\lambda_1 = 0.22$, $\lambda_2 = 0.08$.
We compare, in particular, the second of these values to the prediction interval as estimated by Theorem \ref{the_theorem}, which was calculated to be $(0.18, 0.33)$.  We note especially that the value of $\lambda_1$ is well within the interval  and almost in its center. We also note that the prediction interval does not include the values of $\lambda_0$ and $\lambda_2$.

In the provided example
$\pi_1$ is large compared to $\pi_2$. Both the lower and upper
bound seem sensitive to $\pi_1$, and improve
for $\pi_1$ close to 1. The estimated bounds, we suggest, can be improved by finding
a better estimate for $\|\delta\|$.

\begin{figure}
\begin{tabular}{ccc}
{\includegraphics[width = 0.32\textwidth]{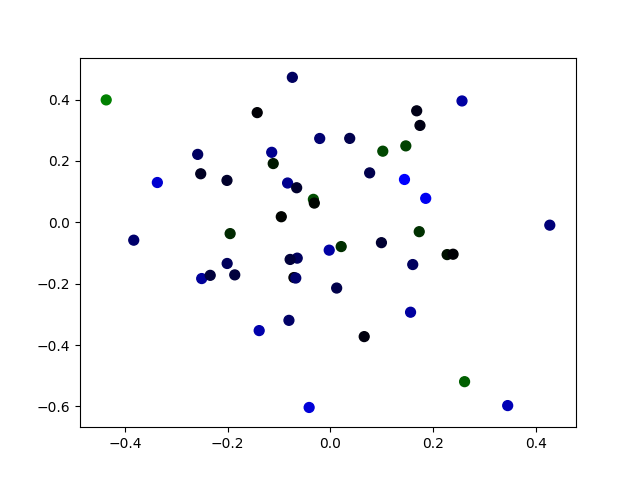}} & 
{\includegraphics[width = 0.32\textwidth]{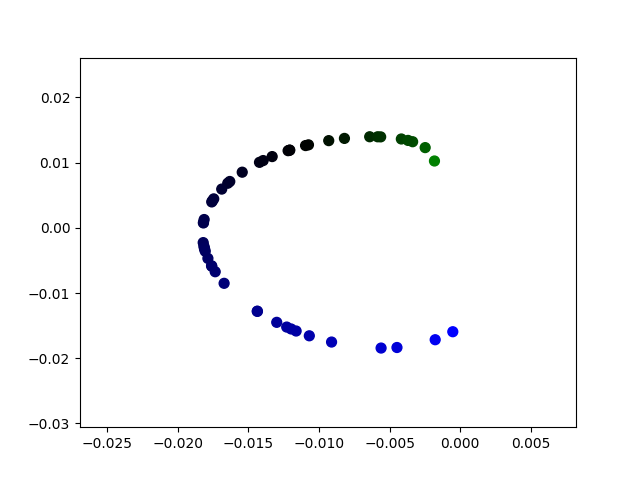}} &
{\includegraphics[width = 0.32\textwidth]{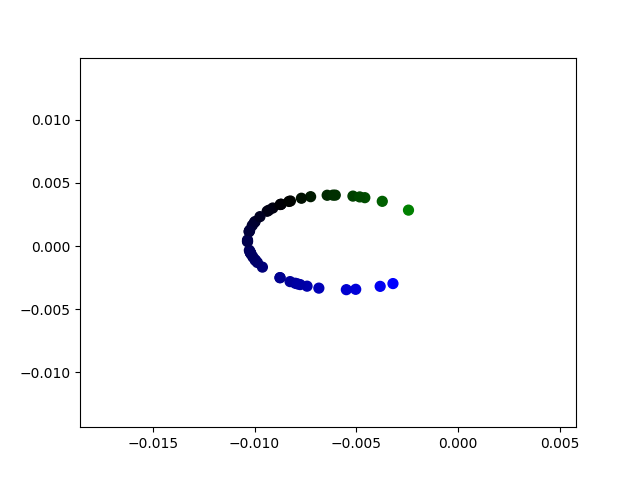}} \\
$a)$ & $b)$ & $c)$ \\
{\includegraphics[width = 0.32\textwidth]{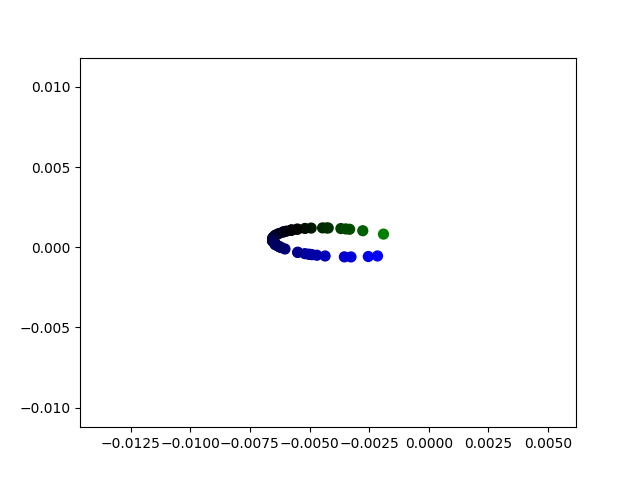}} &
{\includegraphics[width = 0.32\textwidth]{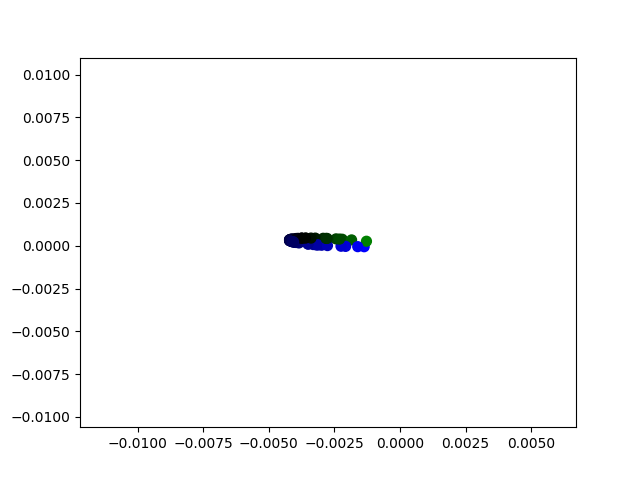}} &
{\includegraphics[width = 0.32\textwidth]{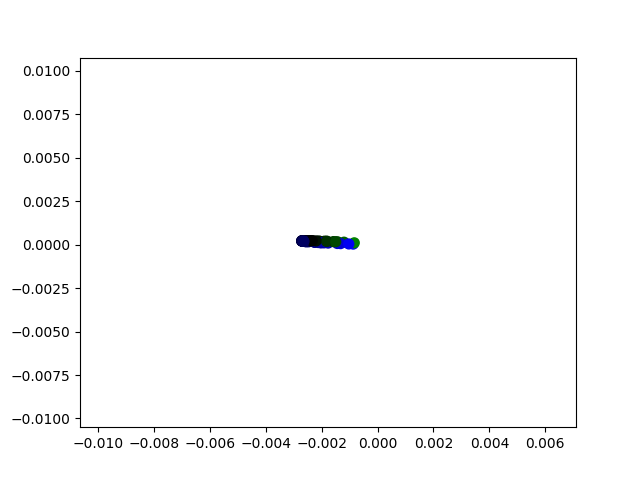}} \\
$d)$ & $e)$ & $f)$ 
\end{tabular}
\caption{Clustering with parameters as in Table \ref{table}.}  Since $\lambda_0$ dominates, the system converges to a cluster before it converges to the origin. \label{seq}
\end{figure}

Fig. \ref{seq} presents clustering with parameters as in Table \ref{table}.
In the sequence of pictures, each frame
is separated from the next by one iteration.

Since the dynamical system $T_x$ is a linear contraction, the vector $y$ will eventually converge to zero. However, Theorem \ref{the_theorem} predicts that the eigenvalue $\lambda_0$ is dominant for the chosen values of parameters,  while all other eigenvalues are quite smaller. This means that long before all points converge to zero, the dynamical system $T_x$ will become two-dimensional (recall, $\cK= K \otimes K$), which experimentally is observed as convergence of all points to a single cluster (see Fig. \ref{seq}). The location of the system in the two-dimensional hyperplane to which it has contracted is, for example, described by the two coordinates of the center of mass of the cluster.

\bgroup
\def\tabcolsep{15pt}
\begin{table}[h!]
\begin{tabular}{| c c c c  c c c c |}
\hline
$n$ & $\pi_1$ & $\pi_2$ & $\mu_1$ & $\mu_2$ & $\sigma_1$ & $\sigma_2$ & $\omega$\\
$50$ & $0.98$ & $0.02$ & $-10$ & $15$ & $1$ & $1$ & $1$\\
\hline
\end{tabular}
\vspace{5pt}
\caption{Simulation parameters for exemplifying Theorem \ref{the_theorem}.}
\label{table}
\end{table}
\egroup

\section{Conclusions and further work}

We have described a particular algorithm which treats a collection of elements of an image as a certain dynamical system.  As input we use a Gaussian mixture distribution of several components. The clustering seems efficient when compared to its overall simplicity. As described, the algorithm preserves, or, respects the underlying data structure.  In general the clustering seems satisfactory given that the deviation of components is not too high as compared to the sizes of the standard deviations. As an adaptation, we suggest a certain redefined model less sensitive to high deviations in the Gaussian mixture distributions. This algorithm produces remarkably efficient clustering even if the components of the mixture distribution are not well isolated (that is the standard deviations are comparable in size to the difference between the means).

We also construct bounds for the second largest eigenvalue of the kernel matrix used in the algorithm for a Gaussian mixture $P = \pi_1 P^1 + \pi_2 P^2$. The bounds seem reasonable given that $\pi_1$ is large. When the choice of the parameters of the Gaussian mixture is such that the upper bound on the second eigenvalue constructed  in Theorem $\ref{second_ev}$ is visibly less then the lower bound on the top eigenvalue, provided in Theorem $\ref{shi}$, then  clustering of such mixture invariably results in a single cluster. We were not able to bound higher eigenvalues rigorously at this point, however, we do observe numerically that in the case when there is a further gap between the second and the third eigenvalues, the dynamical system $T_x(y)$ quickly converges to two-dimensional.

%Both the lower and upper bound seem sensitive to $\pi_1$, and improve for $\pi_1$ close to 1. The estimated bounds, we suggest, can be improved by finding a better estimate for $\|\delta\|$.

\section{Appendix}

\noindent \textbf{Calculating} $\bm r$. We have
\begin{eqnarray}
\nonumber 	r &=& (\pi_1 \pi_2 \iint K(x,y)^2 dP^1(x)dP^2(y))^{1/2} \\
\nonumber &=&\left(\pi_1 \pi_2  { \pi \omega \over \sqrt{4 \sigma_1^2+4 \sigma_2^2+\omega^2 } } e^{-{2 (\mu_1+\mu_2)^2 \over 4 \sigma_1^2+4 \sigma_2^2+\omega^2  }  } \right)^{1 \over 2} \\
\nonumber &=&\sqrt{\pi_1 \pi_2}  { \sqrt{\pi \omega} \over \left(4 \sigma_1^2+4 \sigma_2^2+\omega^2 \right)^{1 \over 4} } e^{-{(\mu_1+\mu_2)^2 \over 4 \sigma_1^2+4 \sigma_2^2+ \omega^2  }  }.
\end{eqnarray}

\medskip	

\noindent {\bf Calculating}  $\bm{\left\| K \right\|_{L^2_{P^2} \times L^2_{P}}}$.
\begin{align*}
  \iint K(x,z)^2 d P^2(x) d P(z)  &= \pi_1 \iint K(x,z)^2 d P^2(x) d P^1(z)  + \\
  & \hspace{20pt} +\pi_2 \iint K(x,z)^2 d P^2(x) d P^2(z)   \\
  &=\pi_1 { \pi \omega \over \sqrt{4 \sigma_1^2+4 \sigma_2^2+\omega^2 }} e^{-{2 (\mu_1+\mu_2)^2 \over 4 \sigma_1^2+4 \sigma_2^2+ \omega^2  }  }+\\
  &\hspace{20pt}+  \pi_2 { \pi \omega \over \sqrt{8 \sigma_2^2+\omega^2 }} e^{-{8 \mu_2^2 \over 8 \sigma_2^2+ \omega^2  }  }.
\end{align*}

\medskip

\noindent \textbf{Calculating} $\bm{\| \phi_1^1 \|_{P^2}}$.  We have
\begin{align*}
	\int \phi_1^1(x)^2 dP^2(x)= \int A_0 e^{-A_1(x-\mu_1)^2}H_1\left(A_2 \frac{x-\mu_1}{\sigma_1}\right)^2 e^{ -A_3(x-\mu_2)^2} d x,
\end{align*}
where 
$$
A_0 =  \frac{(1+2\beta^*)^{1/4}}{2 \sqrt{2\pi \sigma_2^2}}, \quad
A_1  = \frac{\sqrt{1+2\beta^*}-1}{2\sigma_1^2},\quad
A_2  = \left( \frac{1+2\beta^*}{4} \right)^{1/4},\quad
A_3  = \frac{1}{2\sigma_2^2},
$$
%\begin{align*}
%A_0 & =  \frac{1}{\sqrt{2\pi \sigma_2^2}} \left(\frac{(1+2\beta^*)^{1/8}}{\sqrt{2}}\right)^2, \\
%A_1 & = \frac{\sqrt{1+2\beta^*}-1}{2\sigma_1^2},\\
%A_2 & = \left( \frac{1+2\beta^*}{4} \right)^{1/4},\\
%A_3 & = \frac{1}{2\sigma_2^2},
%\end{align*}
and $\beta^* = 2\sigma_1^2/\omega^2$. We calculate
\begin{align*}
&-A_1(x-\mu_1)^2-A_3(x-\mu_2)^2\\
&\hspace{20pt}
= -A_1(x^2-2x\mu_1 + \mu_1^2)
-A_3(x^2-2x\mu_2 + \mu_2^2)\\
&\hspace{20pt} = -(A_1+A_3)
x^2 + (2\mu_1 A_1+2\mu_2 A_3)x -A_1\mu_1^2-A_3\mu_2^2\\
&\hspace{20pt} = 
-(A_1+A_3)\left( x^2 
- 2\frac{ \mu_1A_1+\mu_2 A_3 }{A_1+A_3}x
+\frac{A_1\mu_1^2+A_3\mu_2^2}{A_1+A_3}\right)\\
&\hspace{20pt} = -(A_1+A_3)
\left[ \left(x-\frac{ \mu_1 A_1+\mu_2 A_3 }{A_1+A_3}\right)^2
  - \left( \frac{ \mu_1 A_1+\mu_2 A_3 }{A_1+A_3} \right)^2 + \right. \\
  &\hspace{218pt} \left. + \frac{A_1\mu_1^2+A_3\mu_2^2}{A_1+A_3} \right]\\
&\hspace{20pt} = - a\left[ (x-b)^2 + c \right]
\end{align*}
with $a$, $b$ and $c$ defined appropriately. Thus
\begin{align*}
A_0\int \phi_1^1(x)^2 dP^2(x)&= A_0 \int e^{-a((x-b)^2+c)} H_1\left(A_2 \frac{x-\mu_1}{\sigma_1}\right)^2 dx\\
&= A_0 e^{-a c} \int e^{-a(x-b)^2} H_1\left(A_2 \frac{x-\mu_1}{\sigma_1}\right)^2dx\\
& = A_0 \frac{\sigma_1}{A_2} e^{-ac} \int e^{-a\left( \frac{\sigma_1 y}{A_2}+\mu_1-b\right)^2} H_1^2(y)dy\\
%&= A_0 \frac{\sigma_1}{A_2} e^{-ac} \int e^{-a\left(\frac{\sigma_1}{A_2}\left( y+A_2 \frac{\mu_1-b}{\sigma_1}\right)\right)^2}H_1^2(y)dy\\
%&= A_0 \frac{\sigma_1}{A_2} e^{-ac} \int e^{-a\frac{\sigma_1^2}{A_2^2}\left( y +A_2 \frac{\mu_1-b}{\sigma_1}\right)^2}H_1^2(y)dy\\
&= A_0 \frac{\sigma_1}{A_2} e^{-ac} \int e^{-a\frac{\sigma_1^2}{A_2^2}\left( y
-A_2 \frac{b-\mu_1}{\sigma_1}\right)^2}H_1^2(y)dy.
\end{align*}
Let $\tilde{x} = A_2 (b-\mu_1)/{\sigma_1}$ and $R = a\sigma_1^2/A_2^2$, so the above is
$$
A_0 \frac{\sigma_1}{A_2} e^{-ac} \int e^{-R\left( y -\tilde{x}\right)^2}H_1^2(y)dy.
$$
Since
$$
\int e^{-R\left( y -\tilde{x}\right)^2} H_1^2(y)dy = \frac{ \sqrt{\pi} (2R\tilde{x}^2+1) }{2R^{3/2}},
$$
we get
\begin{align}
	\int (\phi_1^1(x))^2 dP^2(x) =A_0 \frac{\sigma_1}{A_2} e^{-ac} \frac{ \sqrt{\pi} (2R\tilde{x}^2+1) }{2R^{3/2}}. \label{good}
\end{align}
Note that we have made implicit the assumption $R \geq 0$.

\medskip

\textbf{Calculating} $\bm{\lambda_0 \int \phi_1^1(x)\phi_0^1(x) dP(x)}$.
We have 
\begin{align*}
&\lambda_0 \int \phi_1^1(x)\phi_0^1(x) dP(x)\\
&\hspace{20pt} = 
\lambda_0 \int \pi_1 \phi_1^1(x)\phi_0^1(x) dP^1(x)
+ \lambda_0 \int \pi_2 \phi_1^1(x)\phi_0^1(x) dP^2(x)\\
&\hspace{20pt} = 
\lambda_0 \pi_2 \int  \phi_1^1(x)\phi_0^1(x) dP^2(x)\\
&\hspace{20pt} = \int
B_1 e^{-B_2 (x-\mu_1)^2} H_1 \left(B_3(x-\mu_1)\right)  B_4  e^{-B_2 (x-\mu_1)^2} e^{-B_5 (x-\mu_2)^2}dx\\
&\hspace{20pt}
= B_1B_4 \int e^{-2B_2(x-\mu_1)^2-B_5(x-\mu_2)^2}H_1(B_3(x-\mu_1))dx,
\end{align*}
where
\begin{align*}
B_1  = \frac{\lambda_0 \pi_2}{\sqrt{2\pi \sigma_2^2}}  \left(\frac{1+2\beta^*}{\sqrt{2}}\right)^{1/8},\quad B_2  = & \frac{\sqrt{1+2\beta^*}-1}{4\sigma_1^2},\quad B_3  = \left( \frac{1+2\beta^*}{4} \right)^{1/4},\\
B_4  = (1+2\beta^*)^{1/8},&\quad B_5  = \frac{1}{2\sigma_2^2}.
\end{align*}

We have
\begin{align*}
	& -2B_2(x-\mu_1)^2-B_5(x-\mu_2)^2\\
	&\hspace{20pt} = (-2B_2-B_5)x^2 + (2B_2\mu_1+B_5\mu_2)2x
	-2B_2\mu_1^2-B_5\mu_2^2	\\
	&\hspace{20pt} = 
	(-2B_2-B_5)\left( x^2 - 2\frac{ 2B_2\mu_1+B_5\mu_2 }{2B_2+B_5}x
	+ \frac{ 2B_2\mu_1^2+B_5\mu_2^2 }{ 2B_2+B_5 } \right)\\
	&\hspace{20pt}=
	(-2B_2-B_5)\left[ \left( x- \frac{ 2B_2\mu_1+B_5\mu_2 }{2B_2+B_5} \right)^2 - \left( \frac{ 2B_2\mu_1+B_5\mu_2 }{2B_2+B_5} \right)^2 + \right.
	  \\
          &\hspace{228pt}\left. + \frac{ 2B_2\mu_1^2+B_5\mu_2^2 }{ 2B_2+B_5 } \right]\\
	&\hspace{20pt} = -d((x-e)^2+f),
\end{align*}
where $d$, $e$ and $f$ are defined appropriately. Thus
\begin{align*}
	\lambda_0 \int \phi_1^1(x)\phi_0^1(y)dP(x) &= B_1B_4e^{-df}\int e^{-d(x-e)^2} H_1(B_3(x-\mu_1))dx\\
	&= {B_1B_4 \over B_3}e^{-df} 
	\int e^{-d\left( \frac{y}{B_3}+\mu_1-e\right)^2} H_1(y)dy\\
	&= {B_1B_4 \over B_3} e^{-df}
		\int e^{-\frac{d}{B_3^2}\left( y+B_3\mu_1-eB_3\right)^2} H_1(y)dy\\
	&= {B_1B_4 \over B_3} e^{-df}
			\int e^{-\frac{d}{B_3^2}\left( y-\hat x \right)^2}H_1(y)dy
\end{align*}
with $\hat x$ is defined appropriately.
  
Let $S := d/B_3^2$, then since
\begin{align*}
\int e^{-S\left( y-\hat x\right)^2}H_1(y)dy
=\frac{\sqrt{\pi}}{\sqrt{S}} \hat x,
\end{align*}
we get
\begin{align*}
  \lambda_0 \int \phi_1^1(x)\phi_0^1(x) dP(x) &= {B_1 B_4 \over B_3}e^{-df} \frac{\sqrt{\pi}}{\sqrt{S}} \hat x\\
  & = \frac{\lambda_0 \pi_2}{\sqrt{2\pi \sigma_2^2}}\frac{(1+2 \beta^*)^{1 \over 4}}{2^{1 \over 16}}   \frac{ \sqrt{2} }{(1+2\beta^*)^{1 \over 4}} \frac{\sqrt{\pi}}{\sqrt{S}} \hat x  e^{-df}\\
  & = \frac{ \lambda_0 \pi_2}{2^{1 \over 16} \sigma_2}\frac{\hat x}{\sqrt{S}} e^{-df}.
\end{align*}

\medskip

\noindent \textbf{Calculating} $\bm{\|\phi_0^1\|_{P^2}}$. We have
\begin{align*}
\int \phi_0^1(y)^2 dP^2(y) & = \int B_4^2 B_5 e^{-2B_2(x-\mu_1)^2} e^{-B_5(x-\mu_2)^2} dx,
\end{align*} 
where $B_4$ and $B_5$ are defined above.
We concentrate on
\begin{align*}
&-2B_2(x-\mu_1)^2-B_5(x-\mu_2)^2\\
&\hspace{20pt} = (-2B_2 - B_5)x^2 + (2B_2\mu_1 + B_5\mu_2)2x
-(2B_2\mu_1^2 + B_5\mu_2^2)\\
&\hspace{20pt} = (-2B_2-B_5)\left( x^2
- 2\frac{ 2B_2\mu_1 + B_5\mu_2 } { 2B_2+B_5} + 
\frac{ 2B_2\mu_1^2 + B_5\mu_2^2 } { 2B_2+B_5 } \right)\\
&\hspace{20pt}
= (-2B_2-B_5)\left(   \left( x - \frac{ 2B_2\mu_1 + B_5\mu_2 } { 2B_2+B_5} \right)^2 - \left( \frac{ 2B_2\mu_1 + B_5\mu_2 } { 2B_2+B_5} \right)^2 + \right.\\
&\hspace{232pt} \left. +  \frac{ 2B_2\mu_1^2 + B_5\mu_2^2 } { 2B_2+B_5 } \right)\\
&\hspace{20pt} = -u((x-v)^2+w),
\end{align*}
with $u$, $v$ and $w$ defined appropriately. Thus
\begin{align*}
  \|\phi_0^1\|^2_{P^2}  = \int B_4^2 B_5 e^{-u((x-v)^2+w)} d x& = B_4^2 B_5 e^{-uw} \int e^{-u(x-v)^2} d x \\
  &= B_4^2 B_5 e^{-uw} \sqrt{\pi \over u}.
\end{align*}

\medskip

\textbf{Calculating} $\bm{\|F\|^2_{P^1}}$.

\begin{equation*}
\begin{split}
 \|F \|^2_{P^1} & = \int \left( \int \phi_0^1(x) f(x) dP^1(x)\right)^2 \phi_0^1(z)^2 dP^1(z)\\
& = 
\int \left(\int \phi_0^1(x) f(x) d\frac{P(x)-\pi_2 P^2(x)}{\pi_1} -\int \phi_0(x) f(x) dP(x) \right)^2 \phi_0^1(z)^2 dP^1(z)\\
& = 
\int \left[ \int (\phi_0^1(x)  - \phi_0(x)) f(x) dP(x)
+ \left(\frac{1}{\pi_1}-1\right)\int \phi_0^1(x) f(x) dP(x) \right.\\
&\hspace{50pt} \left. - \frac{\pi_2}{\pi_1} \int \phi_0^1(x) f(x) dP^2(x) \right]^2\phi_0^1(z)^2 dP^1(z)\\
& = \int \left( \int (\phi_0^1(x)-\phi_0(x)) f(x) dP(x) + \Delta\right)^2
\phi_0^1(z) dP^1(z) \\
& \leq \| \delta\|^2_{P} \| f\|^2_{P} + 2 \| \delta\|_P \| f\|_P \Delta 
+ \Delta^2\\
& = \| \delta\|^2_{P}  + 2 \| \delta\|_P \Delta 
+ \Delta^2
\end{split}
\end{equation*}
where $\delta(z) = \phi_0^1(z)-\phi_0(z)$, and
\begin{align*}
\Delta & =  \left(\frac{1}{\pi_1}-1\right)\int \phi_0^1(z) f(z) dP(z)
- \frac{\pi_2}{\pi_1} \int \phi_0^1(z) f(z) dP^2(z) \\
& \leq \left(\frac{1}{\pi_1}-1\right) \| \phi_0^1 \|_P \| f\|_P + \frac{\pi_2}{\pi_1} \| \phi_0^1 \|_{P^2} \| f \|_{P^2}\\
& \leq \left(\frac{1}{\pi_1}-1\right) \sqrt{ \pi_1 \| \phi_0^1 \|^2_{P^1} +  \pi_2 \| \phi_0^1 \|^2_{P^2} } + \frac{\sqrt{\pi_2}}{\pi_1} \| \phi_0^1 \|_{P^2} \sqrt{\| f\|^2_P-\pi_1 \| f \|^2_{P^1} }\\
& \leq \left(\frac{1}{\pi_1}-1\right)  \sqrt{ \pi_1 +  \pi_2 \| \phi_0^1 \|^2_{P^2} }  + \frac{\sqrt{\pi_2}}{\pi_1} \| \phi_0^1 \|_{P^2}.
\end{align*}

\medskip

\noindent \textbf{Calculating} $\bm{\epsilon}$. Set
\begin{align*}
\epsilon^2 = \int \left(\pi_2 \int K(x,y) \phi_0^1(y) dP^2(y)\right)^2 dP(x).
\end{align*}
Then
\begin{align*}
	\epsilon^2 \leq \pi_2^2 \|K \|^2_{P^2 \times P} \| \phi^1_0 \|^2_{P^2}.
\end{align*}

\medskip
	
\noindent \textbf{Calculating} $\bm{e}$. We have
\begin{align*}
  |e| & \le  \left| \int \phi_1^1(y)(\phi_0^1(y)-\phi_0(y)) d P(y) \right| +\left| \pi_2 \int \phi_1^1(y) \phi_0^1(y) d P^2(y) \right| \\
  & \leq \|\phi_1^1 \|_P \|\delta  \|_P +\left| \pi_2 \int \phi_1^1(y) \phi_0^1(y) d P^2(y) \right| \\
  & \leq \frac{\epsilon}{t-\epsilon} \sqrt{ \pi_1 + \pi_2 \| \phi_1^1 \|^2_{P^2}} +\left| \pi_2 \int \phi_1^1(y) \phi_0^1(y) d P^2(y) \right|.
   \end{align*}

\bibliographystyle{plain}
\bibliography{biblio}
\end{document}